\documentclass[12pt]{amsart}
\usepackage{times}

\newtheorem{teo}{Theorem}[section]
\newtheorem{lemma}[teo]{Lemma}
\newtheorem{prop}[teo]{Proposition}
\newtheorem{corollario}[teo]{Corollary}
\theoremstyle{definition}

\newtheorem{example}[teo]{Example}
\newcommand{\R}{\mathbb R}

\newcommand{\de}{\partial}
\newcommand{\tr}{\widetilde{\rho}}

\newcommand{\tv}{\widetilde{\varphi}}

\newcommand{\tu}{\tilde{u}}
\newcommand{\tF}{\tilde{F}}
\newcommand{\debar}{\overline{\de}}

\newcommand{\C}{\mathbb C}
\newcommand{\Po}{\mathbb P}
\newcommand{\G}{\mathcal{G}}

\newcommand{\Ko}{\mathcal{K}}

\newcommand{\B}{\mathbb B}
\newcommand{\Sp}{\mathcal Q}
\newcommand{\A}{\mathcal O^{k+\alpha}(\overline{\mathbb D},\C^n)}
\newcommand{\CA}{\mathcal C^{k+\alpha}(\de{\mathbb D},\C^n)}

\newcommand{\Ac}{\mathcal O^{k+\alpha}(\overline{\mathbb D},\C^{n}\times\C^{n-1})}

\newcommand{\Acr}{\mathcal O^{k+\alpha}(\overline{\mathbb D},\R^{n})}

\newcommand{\Aut}{{\sf Aut}(\mathbb D)}
\newcommand{\D}{\mathbb D}
\newcommand{\oD}{\overline{\mathbb D}}
\newcommand{\la}{\langle}
\newcommand{\ra}{\rangle}
\newcommand{\N}{\mathbb N}
\newcommand{\kd}{\kappa_D}
\newcommand{\Hr}{\mathbb H}
\newcommand{\ps}{{\sf Psh}}
\newcommand{\Hess}{{\sf Hess}}
\newcommand{\subh}{{\sf subh}}
\newcommand{\tl}{\tilde{\lambda}}
\newcommand{\gdot}{\stackrel{\cdot}{g}}
\newcommand{\gddot}{\stackrel{\cdot\cdot}{g}}
\newcommand{\fdot}{\stackrel{\cdot}{f}}
\newcommand{\fddot}{\stackrel{\cdot\cdot}{f}}
\def\v{\varphi}
\def\Re{{\sf Re}\,}
\def\Im{{\sf Im}\,}

\theoremstyle{remark}  \newtheorem{nota}[teo]{Remark}
\theoremstyle{remark}  \newtheorem{question}[teo]{Question}
\numberwithin{equation}{section}

\begin{document}
\title[Pluricomplex Poisson kernel]{The pluricomplex Poisson kernel for strongly convex domains}
%tentative title.
\author[F. Bracci, G. Patrizio, S. Trapani]{Filippo Bracci, Giorgio Patrizio, Stefano Trapani}
%\thanks{\rm  Both authors are partially supported by Progetto MIUR di
%Rilevante Interesse Nazionale {\it Propriet\`a geometriche delle
%variet\`a reali e complesse.}}
\address{F. Bracci: Dipartimento di Matematica, Universit\`a di Roma  ``Tor Vergata'', Via
della Ricerca Scientifica 1, 00133 Roma, Italy.}
\email{fbracci@mat.uniroma2.it}
\address{G. Patrizio: Dipartimento di Matematica ``Ulisse Dini'', Universit\`a di Firenze,
Viale Morgagni 67-A, 50134 Firenze, Italy.}
\email{patrizio@math.unifi.it}
\address{S. Trapani: Dipartimento di Matematica, Universit\`a di Roma  ``Tor Vergata'', Via
della Ricerca Scientifica 1, 00133 Roma, Italy.}
\email{trapani@mat.uniroma2.it}
%\subjclass{Primary 32H50, 37F99.}

%\begin{abstract}
%\end{abstract}
\maketitle

\tableofcontents

\section*{Introduction}

In the past decades  the study of pluri-potential theory and of its
applications  played a central role in complex analysis in several
variables. In particular, since the basic work of Siciak \cite{Si}
and Bedford and Taylor \cite{B-T1}, \cite{B-T2} a great effort was
made to understand the complex Monge-Amp\`ere operator and the
associated generalized Dirichlet problems (for instance, see
\cite{De}, \cite{Kl} and references therein).

Let $D\subset \C^n$ be a bounded convex domain with $z_0\in D$. From
the work of Lempert \cite{Le}, \cite{Le3} and Demailly \cite{De} it
turned out that the  following  homogeneous Monge-Amp\`ere equation
\begin{equation}\label{monge-inner-intro}
\begin{cases}
u\in \ps(D)\\
(\de\debar u)^n=0 \quad \hbox{in}\ D\setminus\{z_0\}\\
\lim_{z\to x}u(z)=0\quad \hbox{for all}\ x\in\de D\\
u(z)-\log |z-z_0|=O(1) \quad \hbox{as } z \to z_0
\end{cases}
\end{equation}
has a  solution $L_{D,z_0}$ which is continuous in
$D\setminus\{z_0\}$ (actually it is smooth there if $D$ is strongly
convex with smooth boundary) and unique.

The function $L_{D,z_0}$ shares many properties with the Green
function for the unit disc $\D\subset \C$. For instance, from an
analytic point of view it can be used to reproduce continuous
plurisubharmonic functions (see \cite{De} or Section
\ref{riproduce-sezione}) while from a geometrical point of view, its
level sets are boundaries of Kobayashi balls centered at $z_0$ and
its associated foliation is the singular pencil of complex geodesics
passing through $z_0$ and thus it can be successfully used in
questions such as classification of domains or biholomorphisms (see,
{\sl e.g.}, \cite{Pa1}, \cite{Pa2}, \cite{BDK}). Thus, such function
deserves the name of {\sl pluricomplex Green's function}.

In \cite{B-P} the first and second named authors concentrated in
studying a homogeneous Monge-Amp\`ere equation with a simple
singularity at the boundary. Namely, the following result has been
proved:
\begin{teo}\label{filo-giorgio}
Let $D\subset \C^n$ be a bounded strongly convex domain with smooth
boundary and let $p\in \de D$. The following Monge-Amp\`ere equation
\begin{equation}\label{monge-boundary}
\begin{cases}
u\in \ps(D)\\
(\de\debar u)^n=0 \quad \hbox{in}\ D\\
u<0\quad \hbox{in}\ D\\
u(z)=0 \quad \hbox{for}\ z\in \de D\setminus\{p\}\\
u(z)\approx \|z-p\|^{-1} \quad \hbox{as } z \to p \hbox{
non-tangentially}
\end{cases}
\end{equation}
has a solution $\Omega_{D,p}\in
C^\infty(\overline{D}\setminus\{p\})$ such that
$d(\Omega_{D,p})_z\neq 0$ and $(\de\debar \Omega_{D,p})^{n-1}(z)\neq
0$ for all $z\in \overline{D}\setminus\{p\}$. Moreover the level
sets of $\Omega_{D,p}$ are boundaries of horospheres of $D$ with
center $p$.
\end{teo}

Here $\ps(D)$ denotes the real cone of plurisubharmonic functions in
$D$ and horospheres are the ``limits of Kobayashi balls'' introduced
by Abate \cite{Ab}, \cite{Ab2} and coincide with the sub-level sets
of Busemann functions of  geodesics whose closure contain $p$ (see
\cite{Bus}). The function $\Omega_{D,p}$ has been defined by means
of the boundary spherical representation of Chang-Hu-Lee \cite{CHL}
(see Section \ref{prelimino}). In \cite{B-P}, among other things, it
has been proved that $\Omega_{D,p}$ can be used to characterize
biholomorphisms and that its associated foliation is the fibration
of complex geodesics of $D$ whose closure contain $p$.

The aim of this paper is to  study the  properties of $\Omega_{D,p}$
in depth. We will show that $\Omega_{D,p}$ shares many properties
with the Poisson kernel for the unit disc $\D$ and therefore it
deserves the name of {\sl pluricomplex Poisson kernel} of $D$ with
singularity at $p\in \de D$.

More in detail, we show that a version of the classical
Phragmen-Lindel\"of theorem on the growth of subharmonic functions
in $\D$ holds for plurisubharmonic functions in $D$, proving that
$\Omega_{D,p}$ is the maximal element of the family
\begin{equation*}
\begin{cases}
u \in \ps (D) \\
\limsup_{z\to x} u(z)\leq  0 \quad \hbox{for all}\ x\in \de D\setminus\{p\}\\
\displaystyle{\liminf_{t\to 1} |u(\gamma(t))(1-t)|\geq 2\Re(\la
\gamma'(1),\nu_p\ra^{-1})},
\end{cases}
\end{equation*}
where $\nu_p$ is the unit outward normal to $\de D$ at $p$ and
$\gamma$ is any $C^1$-curve in $D$ such that $\gamma(1)=p$ and
$\gamma'(1)\not\in T_p\de D$ (see Section \ref{estremalita}). In due
course we will find the exact behavior of $\Omega_{D,p}(z)$ as $z$
goes to $p$ along  non tangential directions to $\de D$ at $p$ (see
Corollary \ref{oicomevobene}).

Next, we deal with uniqueness properties of $\Omega_{D,p}$. These
are essentially of two types: analytic and geometric. From an
analytic point of view we show that $\Omega_{D,p}$ is the only
solution of the homogeneous Monge-Amp\`ere equation which is zero on
$\de D\setminus\{p\}$ and behaves like $\Omega_{D,p}$ as $z$ tends
to $p$ (see Theorem \ref{unico-comportamento-al-bordo}). This is the
analogous of the uniqueness statement for the pluricomplex Green
function, except that the behavior of $\Omega_{D,p}$ near $p$ is
universal only along non tangential directions, but it might depend
on the domain $D$ itself along  other directions. From a geometrical
point of view we show that $\Omega_{D,p}$ is the only $C^2$ solution
(up to multiplication by constants) of the homogeneous
Monge-Amp\`ere equation which is zero on $\de D\setminus\{p\}$ and
whose associated foliation is the fibration of $D$ in complex
geodesics whose closure contain $p$ (see Theorem
\ref{foliazione-unica}). This fact is then used to show a couple of
interesting other characterizations of $\Omega_{D,p}$ both in terms
of its level sets (see Proposition \ref{level-unico}) and in terms
of its behavior under pull-back with holomorphic self-maps of $D$
(see Proposition \ref{contra-unico}).

We also show in Theorem \ref{GvP} that $L_{D,z_0}$ and
$\Omega_{D,p}$ have the same relationship as the Green function and
the Poisson kernel in $\D$, namely
\begin{equation}\label{relaZ}
\Omega_{D,p}(z_0)=-\frac{\de L_{D,z_0}}{\de \nu_p}(p).
\end{equation}
This  is used to write explicitly the ``noyaux de Poisson
pluricomplexes canonique'' of Demailly \cite{De} and, applying his
theory, to obtain a somewhat explicit reproducing formula for
continuous plurisubharmonic functions of $D$ in terms of $L_{D,z_0}$
and $\Omega_{D,p}$ (see Theorem \ref{riproduce}). In particular, for
pluriharmonic functions $F\in C^0(\overline{D})$ we obtain the
following formula which is the analogous of that for harmonic
functions in the disc:
\[
F(z)=\int_{p\in \de D} |\Omega_{D,p}(z)|^n F(p) \omega_{\de D}(p),
\]
where $\omega_{\de D}$ is a positive real $(2n-1)$-form on $\de D$
which depends only on $D$.

As a spin off result, using the properties of $\Omega_{D,p}$, we
also prove that horospheres are (smooth and) strongly convex away
from their center (see Theorem \ref{strettaconv}).

The proofs of the previous properties of $\Omega_{D,p}$ are based on
a mix of different techniques. In particular we will make a strong
use of families of complex geodesics  and their regularity
properties. Thus in Section \ref{regulo-sezione} we deal with
regularity for such families gathering some known but disperse
information and proving the precise results needed for our
arguments. In particular, using a suitable ``attached analytic
discs'' approach, we prove (Theorem \ref{regolaritaG}) that the set
of complex geodesics in $D$  is a finite dimensional closed
submanifold in the open set of the complex Banach space $\A$  made
of non-constant holomorphic attached discs whose first $k$-th
derivatives extend $\alpha$-H\"older continuous on $\de \D$. This
result, interesting on its own, allows to obtain stability and
regularities properties for families of complex geodesics (Section
\ref{regulo-sezione}) and for their {\sl Lempert's projections},
that is, the holomorphic retractions of $D$ with affine fibers onto
complex geodesics introduced by Lempert in \cite{Le} which will play
a fundamental role in our discussion (see Section
\ref{retratto-sezione}).

The plan of the paper is as follows. In the first section we recall
some preliminaries about complex geodesics, the boundary spherical
representation of Chang, Hu and Lee \cite{CHL} and the results in
\cite{B-P} as needed to make this work as self-contained as
possible. In section two we deal with regularity for families of
complex geodesics by studying their differential properties  and, as
a corollary of our construction, we recover with a different proof
some stability results by Huang \cite{Hu}, \cite{Hu2}. In the third
section we study Lempert's projections. We first show that
holomorphic retractions on a given complex geodesics are not unique
but Lempert's projections can be characterize as  the unique
retractions with affine fibers. Then we examine the variation of
Lempert's projections with respect to boundary data and prove
regularity.  In section four we investigate the shape of
horospheres. We prove that they are strongly convex away their
center (where they are $C^{1,1}$) using Jacobi vector fields. In the
fifth section we state and prove the Phragmen-Lindel\"of theorem for
strongly convex domains and we compute the limits of $\Omega_{D,p}$
along non complex-tangential directions. In the sixth section we
prove \eqref{relaZ} and in section seven we deal with uniqueness.
Finally, in section eight we recall Demailly's theory for
reproducing plurisubharmonic functions and find the explicit
reproducing formulas using $\Omega_{D,p}$.

\section{Preliminaries}\label{prelimino}

Let $D$ be a bounded strongly convex domain in $\C^n$ with smooth
boundary.  A {\sl complex geodesic} is a holomorphic map $\v:\D\to
D$ which is an isometry between the Poincar\'e metric of
$\D=\{\zeta\in \C: |\zeta|<1\}$ and the Kobayashi distance $k_D$ in
$D$.

According to  Lempert (see \cite{Le} and \cite{Ab}), any complex
geodesic extends smoothly to the boundary of the disc and $\v(\de
\D)\subset \de D$. Moreover,  given any two points $z,w\in
\overline{D}$, $z\neq w$, there exists a complex geodesic $\v:\D\to
D$ such that $z,w\in \v(\oD)$. Such a geodesic is unique up to
pre-composition with automorphisms of $\D$. Also, if $z\in
\overline{D}$ and $v\in \C^n\setminus\{O\}$ (and $v\not\in T_z\de D$
if $z\in \de D$) there exists a unique (still, up to pre-composition
with automorphisms of $\D$) complex geodesic $\v:\D\to D$ such that
$z\in \v(\oD)$ and $\v(\oD)$ is parallel to $v$ (in case $z, w \in
\de D$ this follows from Abate~\cite{Ab4} and Chang, Hu and
Lee~\cite{CHL}). In case $z\in D$ and $w\in \overline{D}$, $w\neq
z$, ({\sl respectively} $v\in T_zD$) one can   choose uniquely a
complex geodesic $\v:\D\to D$ requiring that $\v(0)=z$ and $\v(t)=w$
for some $0<t\leq 1$, with $t=1$ if and only if $w\in \de D$ ({\sl
respect.} $\v'(0)=t v$ for some $t>0$). With an abuse of notation,
when no risk of confusion arises, we call ``complex geodesic'' also
the image of a complex geodesic $\varphi:\D\to D$.

If $\v:\D\to D$ is a complex geodesic then there exists a
holomorphic map $\tv: \D \to \C^n$, called the {\sl dual map} of
$\v$, such that $\tv$ extends smoothly to $\de \D$ and
$\tv(e^{i\theta})=e^{i\theta}\mu(e^{i\theta})\de
r_{\v(e^{i\theta})}$, with $r$ being a defining function of $D$ near
$\varphi(\de \D)$ and $\mu>0$ normalized so that
\begin{equation}\label{norm-1}
\tv(\zeta)\cdot\v'(\zeta)\equiv 1
\end{equation}
for all $\zeta\in \D$ (see \cite{Le}).

Let $\v:\D\to D$ be a complex geodesic. In \cite{Le1} and \cite{Le2}
(see also Pang \cite{Pg1}) Lempert defines a biholomorphic change of
coordinates $G:D\to D'$ which ``linearizes'' $\v$. Namely, he proves
that  $G$ extends smoothly on $\de D$, that $G\circ
\v(\zeta)=(\zeta,0,\ldots,0)$ and $\widetilde{G\circ
\v}(\zeta)=(1,0,\ldots,0)$. The domain $D'=G(D)$ is no longer convex
in general but it is strictly linearly convex near $G(\v(\de \D))$,
namely, the real Hessian of any defining function of $D'$ is
positive on the complex tangent space at any point of $\de D'$  near
$G(\v(\de \D))$. In the rest of the paper we will refer to such a
$G$ as  the Lempert biholomorphism which linearizes $\v$.

Considering the foliation of all  complex geodesics passing through
a given point $z_0\in D$, Lempert constructed a map $\Phi_{z_0}:D\to
\B^n$, called {\sl spherical representation} of $D$ at $z_0$, which
is defined by $\Phi_{z_0}(z)=\zeta \v_z'(0)/\|\v_z'(0)\|\in \B^n$
where $\v_z:\D\to D$ is a complex geodesic such that $\v_z(0)=z_0$,
$\v_z(\zeta)=z$ for $z\neq z_0$ and $\Phi_{z_0}(z_0)=O$. The map
$\Phi_{z_0}$ which is continuous on $D$, extends $C^\infty$ on
$\overline{D}\setminus\{z_0\}$. In his work \cite{Le} Lempert proved
that $L_{D,z_0}:=\log\|\Phi_{z_0}\|$ solves
\eqref{monge-inner-intro}.

Similarly, considering all complex geodesics whose closure contain a
given boundary point $p\in \de D$, Chang, Hu and Lee (see
\cite[Theorem~3]{CHL}) constructed a {\sl boundary spherical
representation}. For the reader convenience and since it will be
useful later, we recall here the construction of Chang, Hu and Lee
 as needed for our aim. Let $p\in \de D$ and let $\nu_p$
be the unit outward normal to $\de D$ at $p$.  Denote
$$L_p:=\{v \in \C^{n} |  \|v\|=1,\ \la v, \nu_p\ra >0,\ iv\in T_p\de D \}$$
and let $v\in L_p$. In what follows we will say that a complex
geodesic $\varphi_v: \D \to D$ whose closure contains the point
$p\in \de D$ is in {\sl the Chang-Hu-Lee normal parametrization}
(with respect to $v\in L_p$) if $\varphi(1)=p$ and $\varphi'(1)=\la
v, \overline{\nu_p}\ra v$ and $\Im \la \v''(1),\overline{\nu_p}\ra
=0$. In~\cite{CHL} Chang, Hu and Lee proved that for all $v\in L_p$
there exists a {\sl unique} complex geodesic in {\sl the
Chang-Hu-Lee normal parametrization} with respect to $v$.

Up to rigid movements of $\C^n$,   assume  that
$\nu_p=e_1=(1,0,\ldots, 0)$ and thus $L_p$ reduces to
$L_p=\{v=(v_1,\ldots, v_n)\in \C^n: \|v\|=1, v_1>0\}$. For any $v\in
L_p$ the map $\eta_v:\D \ni \zeta \mapsto e_1+(\zeta-1)v_1 v$ is a
complex geodesic of $\B^n$, $\eta_v(1)=e_1$ and $\eta_v'(1)=v_1v$.
Then the {\sl boundary spherical representation} $\Phi_p: D\to \B^n$
is defined as follows:
\[
\Phi_p(z)=e_1+(\zeta_z-1)v_{z,1}v_z,
\]
where $\zeta_z\in \D$ and $v_z\in L_p$ are the unique data such that
$\varphi_{v_z}(\zeta_z)=z$. The map $\Phi_p$ is a smooth
diffeomorphism whose inverse is
$\Phi_p^{-1}(w)=\varphi_{v_w}(\zeta_w)$, where $\zeta_w \in \D$ and
$v_w\in L_p$ are the unique data such that $w=\eta_{v_w}(\zeta_w)$.
Moreover $\Phi_p, \Phi_p^{-1}$ extend continuously up to the
boundary and $\Phi_p(p)=e_1$. In particular it follows that $\Phi_p$
is holomorphic on all complex geodesics in $D$ whose closure contain
$p$ and sends such complex geodesics to complex geodesics in $\B^n$
whose closure contain~$e_1$.

Following Abate (\cite{Ab}, \cite{Ab2}) we define a {\sl horosphere}
$E_D(p,z_0, R)$ of center $p\in \de D$, pole $z_0\in D$ and radius
$R>0$ as
\[
E_D(p,z_0, R):=\{z\in D: \lim_{w\to
p}[k_D(z,w)-k_D(z_0,w)]<\frac{1}{2}\log R\}.
\]
The limit in the definition of $E_D(p,z_0, R)$ exists since $D$ is
strongly convex and  any such horosphere $E_D(p,z_0, R)$ is a
sub-level set of the Busemann function of any geodesic whose closure
contains $p$ (see \cite{Trap}).

In \cite[Corollary 6.2]{B-P} it was proved that {\sl $\Phi_p$ maps
horospheres of $D$ centered at $p$ onto horospheres of $\B^n$
centered at $e_1$}, which, since horospheres of $\B^n$ are complex
ellipsoid, implies  in particular that the boundaries of horospheres
are smooth away from the center $p$.

Let $\Omega_{\B^n, e_1}(z)=-\frac{1-\|z\|^2}{|1-z_1|^2}$. The
sub-level sets of $\Omega_{\B^n, e_1}$ corresponds to horospheres of
$\B^n$ with center $e_1$ and pole $O$ (see, {\sl e.g.}, \cite{Ab},
\cite{Ab2}). In \cite{B-P} we defined
\[
\Omega_{D,p}:=\Omega_{\B^n, e_1}\circ \Phi_p
\]
and proved Theorem \ref{filo-giorgio}. For further use we notice
that
\[
E_D(p, \Phi_p^{-1}(O), R)=\{z\in D: \Omega_{D,p}(z)<-\frac{1}{R}\}.
\]
Finally,  let
\[
P(\zeta):=\frac{1-|\zeta|^2}{|1-\zeta|^2}
\]
be the Poisson kernel on $\D=\{|\zeta|<1\}$. Recall that $P$ is
harmonic in $\D$,  $\lim_{\zeta\to x}P(\zeta)=0$ for $x\in \de
\D\setminus\{1\}$ and $\lim_{\R\ni r\to 1^-}P(r)(1-r)=2$.
 From the very
definition it follows that for all $v\in L_p$
\begin{equation}\label{comefo}
\Omega_{D,p}(\v_v(\zeta))=-P(\zeta)/v_1^2.
\end{equation}

\section{Regularity for families of complex
geodesics}\label{regulo-sezione}

In this section we  state some results about regularity of families
of complex geodesics in strongly convex domains which we need later.
From these we also rediscovered some facts already known or
implicitly contained in other papers such as \cite{Le}, \cite{Le1},
\cite{Hu}, \cite{Hu2}). Our presentation owes much to the works
\cite{Glo}, \cite{Tre}, \cite{Trap2}.

In all this section $D$ will be a bounded strongly convex domain of
$\C^n$ with smooth boundary. Given  $k\geq 2$ and $\alpha\in (1/2,
1)$ we denote by $\A$ the set of all holomorphic maps from $\D$ to
$\C^N$ which extends $C^k$ on $\overline{\D}$ and such that their
$k$-th derivatives are $\alpha$-H\"older on $\overline{\D}$ (a map
$f:\oD\to \C^n$ is $\alpha$-H\"older if there exists $C>0$ such that
$\|f(\zeta_0)-f(\zeta_1)\|\leq C|\zeta_0-\zeta_1|^\alpha$ for all
$\zeta_0, \zeta_1\in \D$). The set $\A$ is a complex Banach space
when endowed with the  norm
\[
\|f\|_{k+\alpha}=\sum_{j=1}^k
\sup_{\zeta\in\de\D}\|f^{(j)}(\zeta)\|+\sup_{\zeta_0, \zeta_1\in
\oD, \zeta_0\neq \zeta_1}
\frac{\|f^{(k)}(\zeta_0)-f^{(k)}(\zeta_1)\|}{|\zeta_0-\zeta_1|^\alpha}.
\]
Let $\G$ be the set of complex geodesics from $\D$ to $D$. By
Lempert's theory \cite{Le} it follows that $\G\subset \A$. Let also
denote by $M\subset \A$ the set of constants with value in $\de D$.
It is clear that $M$ is a closed set in $\A$.

\begin{teo}\label{regolaritaG}
The set $\G$ is a closed submanifold of $\A\setminus M$ of real
dimension $4n-1$.
\end{teo}
\begin{proof}
Let $\{f_n\}\subset \G$ and assume that $f_n\to f$ in $\A$. Since
the domain $D$ is strongly (pseudo)convex then either $f(\D)\subset
D$---and from the continuity of $k_D$ it follows easily that $f\in
\G$ as well---or $f\in M$. Thus $\G$ is closed in $\A\setminus M$.
Let $f_0\in \G$. We want to prove that $\G$ is a submanifold of $\A$
near $f_0$.

Let $G:D\to D'=G(D)$ be the Lempert biholomorphisms which linearizes
$f_0$. Then $G\circ f_0(\zeta)=(\zeta,0,\ldots,0)$ and the dual map
$\widetilde{G\circ f}_0(\zeta)\equiv (1,0,\ldots, 0)$. Notice that
$G$ extends $C^\infty$ up to $\de D$. Thus we can extend
(arbitrarily) $G|_{\de D}$ to some $C^\infty$ map, denoted by
$\tilde{G}$, from $\C^n$ to $\C^n$. We have thus a morphism $\Lambda
: \CA\to \CA$ given by $\Lambda (f)=\tilde{G} \circ f$. The morphism
$\Lambda$ is $C^\infty$ and maps the set of complex geodesics of $D$
onto the set of complex geodesics of $D'$. Assume for the moment
that we proved that $\Lambda(\G)$ is a finite dimensional
submanifold of $\A$ near $G\circ f_0$, and thus a finite dimensional
submanifold of $\CA$. Repeating the argument with $G^{-1}$, we find
a $C^\infty$ map $\Lambda':\CA\to \CA$ such that $\Lambda \circ
\Lambda'|_{\Lambda(\G)}={\sf Id}|_{\Lambda(\G)}$. Thus
$\Lambda'|_{\Lambda(\G)}$ is an embedding with
$d\Lambda'(T\Lambda(\G))$ finite dimensional, thus closed and
complemented in $\CA$. Therefore $\G=\Lambda'(\Lambda(\G))$ is a
finite dimensional submanifold (see, {\sl e.g.}, \cite{AMR}). We are
then left to show that $\Lambda(\G)$ is a finite dimensional
submanifold.

Thus, we can assume from the beginning that
$f_0(\zeta)=(\zeta,0,\ldots,0)$ and $\tilde{f}_0(\zeta)=(1,0,\ldots,
0)$ in $D$---here however the domain $D$ is no longer strongly
convex, but it is   strictly linearly convex near $f_0(\de \D)$. By
the very definition of the dual map and by \eqref{norm-1} it follows
that if $f$ is a complex geodesic of $D$ close to $f_0$ in $\A$ then
$\tilde{f}$ is close to $\tilde{f_0}$ in $\A$, where, with some
abuse of notation, we identify the one form $\tilde{f}$ with the
vector of its components.

Let $\Po^{n-1}(\C)$ be the space of complex hyperplanes passing
through the origin $O$. Let $\Psi: \de D \to \C^n \times
\Po^{n-1}(\C)$ be defined by $\Psi(p):=(p, T^{\mathbb C}_p \de D)$.
Let $S=\Psi (\de D)$. By the very definition $(f_0,
[\tilde{f_0}])(\de \D)\subset S$. Moreover, since $\de D$ is
strongly pseudoconvex near $f_0(\de D)$, then $S$ is a compact
maximal totally real submanifold of $\C^n\times \Po^{n-1}(\C)$ near
$\Psi(f_0(\de D))$ (see \cite{We}).

Let $(z_1,\ldots, z_n)$ be coordinates in $\C^n$ and let
$[z_1:\ldots :z_n]$ be the corresponding homogeneous coordinates in
$\Po^{n-1}(\C)$; that is, the point $[z_1:\ldots :z_n]$ corresponds
to the hyperplane $\{ v=(v_1,\ldots, v_n)\in \C^n : \sum_{j=1}^n
v_j\cdot z_j=0\}$. Let $U_1:=\{[z]\in \Po^{n-1}(\C): z_1\neq 0\}$ be
the chart obtained by identifying $\C^{n-1}$ with $U_1$ via
$(w_1,\ldots, w_{n-1})\to [1:w_1:\ldots: w_{n-1}]$ and let
$R:\C^n\times \C^{n-1}\to \R^{2n-1}$ be a defining function for
$S\cap \C^n\times U_1$ (such a defining function can be easily
defined starting from a global defining function of $D$ in $\C^n$).
Let us consider $\Sp=\{F=(f,g)\in \Ac : R(f,g)|_{\de \D}\equiv 0\}$.
In other words, $F=(f,g)\in \Sp$ if and only if $(f, [1:g])(\de
\D)\subset S$. In particular $(f_0, 0)\in \Sp$. Note that if $f$ is
a complex geodesic close to $f_0$ with dual map
$\tilde{f}=(\tilde{f}_1, \tilde{f}_2)\in \C\times \C^{n-1}$ then
$\min_{\zeta\in \oD}|\tilde{f}_1(\zeta)|>0$ and therefore $(f,
\tilde{f}_2/\tilde{f}_1)\in \Sp$. Conversely, if $(f,g)\in\Sp$ and
$(f,g)$ is close to $(f_0,0)$ in $\Ac$ then
$|f_1'(\zeta)+\sum_{j=2}^n f_j'(\zeta)g_j(\zeta)|>0$ for all
$\zeta\in \oD$ and then $f$ is a stationary disc in $D$ with dual
map $(1,g)/(f_1'+\sum_{j=2}^n f_j'g_j)$ that is, a complex geodesic.
It should be remarked that in this argument one cannot refer
directly to Lempert's theory  because $D$ is not strongly convex in
general. However, since $\de D$ is strongly pseudoconvex near
$f_0(\de \D)$ then for  $f$ close to $f_0$ in $\A$, one can use
Pang's results to relate stationarity to extremality, see
\cite[Section~2]{Pg1}.

The previous discussion shows that there exists a open neighborhood
$W_0\subset \Ac$ of $(f_0,0)$ such that $\pi_1:\Sp\cap
W_0\to\G\cap\pi_1(W_0)$ is bijective, where $\pi_1$ is the
projection on the first factor, namely $\pi_1(f,g):=f$. The map
$\pi_1|_{\Sp\cap W_0}$ is clearly $C^\infty$, and its inverse is
$C^\infty$ as well, being given by $f\mapsto
(f,\tilde{f}_2/\tilde{f}_1)$ with
$\tilde{f}=(\tilde{f}_1,\tilde{f}_2)\in \C\times\C^{n-1}$ the dual
of $f$. Thus $d\pi_1|_{T\Sp}$ is injective and its image is finite
dimensional and hence closed and complemented in $\A$. Therefore, if
we prove that $\Sp$ is a finite dimensional submanifold of $\Ac$
near $(f_0,0)$ then the claim on $\G$ will follow.
 To prove that $\Sp$ is a submanifold  by means of the
implicit function theorem in Banach spaces, it is enough to show
that $dR_{f_0}:\Ac \to \Acr$ is surjective and its kernel
complemented in $\Ac$. We have $dR_{f_0}(f(\zeta))=2 \Re A f(\zeta)$
with $A$ being the $(2n-1)\times (2n-1)$ matrix with entries
$\frac{\de R_j}{\de z_k}(f_0)$. Since $S$ is maximal totally real,
arguing  as in \cite[Theorem 3.1, Lemma 3.2]{S-T} one can prove that
all the Birkhoff partial indices of the operator $f\mapsto 2 \Re A
f(\zeta)$ are  $\geq 1$ and thus, by \cite{Glo} (see also \cite{Tre}
and \cite{Trap2}) $dR_{f_0}$ is surjective. Notice that the
computation of Birkhoff partial indices in \cite[Lemma 3.2]{S-T} was
proved under the assumption that $\de D$ is strongly convex. It is
easy to check that in fact such result holds for strictly linearly
convex domains and therefore it can be used here as, in Lempert's
coordinates, the domain is strictly linearly convex near $f_0(\de
\D)$. Finally, a direct computation (or see \cite{Trap}) shows that
its kernel has finite (real) dimension $4n-1$ and therefore $\Sp$ is
a submanifold of dimension $4n-1$ near $f_0$.
\end{proof}

Let $\kd$ be the Kobayashi metric in $D$. According to Lempert
(\cite{Le}, \cite{Le2}) the map $D\times (\C^{n}\setminus\{O\})\to
\R$ given by  $(z,v)\mapsto \kd(z;v)$ is $C^\infty$. Moreover, since
$\kd(z,\lambda v)=\lambda \kd(z,v)$ for all $(z,v)\in D\times
(\C^{n}\setminus\{O\})$ and $\lambda>0$ it follows that
$d(\kd)_{(z,v)}\neq 0$ for all $(z,v)\in D\times
(\C^{n}\setminus\{O\})$. Therefore the set
\[
\Ko=\{(z,v)\in D\times
(\C^{n}\setminus\{O\}): \kd(z;v)=1\}
\]
 is a $(4n-1)$-real dimensional
submanifold of $D\times (\C^{n}\setminus\{O\})$.

\begin{teo}\label{quala}
The map $V:\G\to \Ko$ defined by $V:f\mapsto (f(0), f'(0))$ is a
diffeomorphism.
\end{teo}
\begin{proof}
By the uniqueness of complex geodesics \cite{Le}, the map $V$ is
bijective. Since $V$ is the restriction of a linear bounded map from
$\A$ to $\C^{2n}$ then it is linear and $C^\infty$. By \cite[Theorem
5]{Le2} the inverse $V^{-1}$ is $C^\infty$ as well and hence $V$ is
a diffeomorphism.
\end{proof}

From this result we obtain some corollaries which will be useful
later on.

\begin{corollario}\label{appendeprop}
Let $\{f_n\}\subset \G$ be such that $f_n\to f$ uniformly on
compacta of $\D$. If $f$ is not constant then $f\in \G$ and
$f_n^{(j)}\to f^{(j)}$ uniformly on $\oD$ for all $j=0,1,\ldots$.
\end{corollario}
\begin{proof}
By Theorem \ref{regolaritaG} if $f$ is not constant then it belongs
to $\G$. Thus, $f_n\to f$ uniformly on compacta of $\D$ implies that
$f_n(0)\to f(0)$ and $f_n'(0)\to f'(0)$. By Theorem \ref{quala}it
follows that $f_n\to f$ in $\A$ for all fixed $k\in \N$. In
particular $f_n^{(j)}\to f^{(j)}$ for all $j=0,1,\ldots$.
\end{proof}

\begin{corollario}\label{appende}
It $(0,1)\ni t\mapsto f_t\in \G$ is a family of complex geodesics
such that $t\mapsto f_t(0)$ and $t\mapsto f'_t(0)$ are $C^\infty$
then $t\mapsto f_t$ is $C^\infty$ in $\A$. In particular the map
$\zeta\mapsto \frac{\de^j f_t}{\de t^j}(\zeta)$ is smooth on $\oD$
for all $j=1,2,\ldots$
\end{corollario}

\begin{lemma}\label{propG}
The map $\G \ni f\mapsto f(0)\in D$ is proper.
\end{lemma}
\begin{proof}
If $\{z_n\}\subset D$ is such that $z_n\to z\in D$ let $f_n\in \G$
be such that $f_n(0)=z_n$. Let  $\{f_{n_k}\}$ be a converging
subsequence. Since the Kobayashi distance is continuous on $D$ it
follows that the limit $f$ of $\{f_{n_k}\}$ is not constant. Then by
Corollary \ref{appendeprop}  it follows that $f\in \G$ and $f(0)=z$.
Hence the map $f\mapsto f(0)$ is proper.
\end{proof}

As a straightforward corollary of Lemma \ref{propG} we have the
following result, first proved with different methods by Huang
\cite[Proposition 1]{Hu2}:

\begin{prop}\label{huang}
Let $c>0$ and let $\G_c:=\{f\in \G: \hbox{dist}(f(0), \de D)\geq
c\}$. Then there exists $c'>0$ such that $\|f\|_{k+\alpha}\leq c'$.
\end{prop}

\section{Lempert's
projections}\label{retratto-sezione}

Let $D\subset \C^n$ be a bounded strongly convex domain with smooth
boundary and let $\v:\D\to D$ be a complex geodesic. According to
Lempert (\cite{Le}, \cite{Le1}, \cite{Le2}), for all $z\in D$ the
equation $\tv(\zeta)\cdot(z-\v(\zeta))\equiv 0$ in the unknown
$\zeta\in \D$  has a unique solution $\zeta:=\tilde{\rho}(z)$. The
map $\tilde{\rho}:D\to \D$ is holomorphic, extends smoothly on $\de
D$ and it is called the {\em left inverse } of $\varphi$ for it
satisfies $\tilde{\rho} \circ \varphi = {\sf id}_{\D}$. By the very
definition
\begin{equation}\label{ro-duale}
\tv(\tilde{\rho}(z))\cdot(z-\v(\tilde{\rho}(z)))\equiv 0.
\end{equation}

\begin{nota}\label{zero-ro}
Let $z\in \overline{D}$. If $\zeta\in \oD$ is such that
$\tv(\zeta)\cdot(z-\v(\zeta))=0$ then $\tr(z)=\zeta$. Indeed, by the
strong convexity of $\de D$, if $z\in \overline{D}\setminus\v(\de
\D)$ then the winding number of the function $\de \D \ni \zeta
\mapsto \tv(\zeta)\cdot(z-\v(\zeta))$ is $1$ (see \cite{Le1},
\cite{Le2}) hence $\zeta=\tr(z)$. On the other hand, if
$z=\v(e^{it})$ for some $t\in \R$, by continuity of $\tr$ it follows
that $\tr(\v(e^{it}))=e^{it}$.  Suppose by contradiction that
$\tv(\zeta)\cdot(\v(e^{it})-\v(\zeta))=0$ for some $\zeta\in
\oD\setminus\{e^{it}\}$. Since the domain is strongly convex the
interior of the real segment $\ell$ joining $\v(e^{it})$ to
$\v(\zeta)$ is contained in $D$.  Then the segment $\ell$
 belongs to the fiber of $\tr$ at $\v(\zeta)$ and, since $\tr$ is
continuous on $\overline{D}$, it follows that
$\tr(\v(e^{it}))=\zeta$ which contradicts $\tr(\v(e^{it}))=e^{it}$.
\end{nota}

Let $\v$ be a complex geodesic and let $\tilde{\rho}$ be its
left-inverse. The map $\rho:D\to \varphi(\D)\subset D$ defined as
$\rho:=\v\circ\tilde{\rho}$ is a holomorphic retraction on $\v(\D)$,
{\em i.e.}, $\rho$ is a holomorphic self-map of $D$ such that $\rho
\circ \rho =\rho$ and $\rho(z)=z$ for any $z \in\varphi(\D)$. It
extends smoothly to $\de D$ and it is called the {\em Lempert
projection} associated to $\varphi$.  The triple $(\varphi, \rho,
\tilde{\rho})$ is the so-called {\em Lempert projection device}. As
remarked for instance  in~\cite[p. 145]{Br}  the Lempert projection
$\rho$  depends only on the image $\varphi(\D)$.

In this section we study regularity of Lempert's left-inverse.
Before that, we make some comments on holomorphic retractions on
strongly convex domains. We start with an example which shows that
there exist infinitely many holomorphic retractions:

\begin{example}
Let $f_{jk}:\B^n\to \D$ be holomorphic functions, $j,k=2,\ldots, n$
and let $\epsilon<1/2n$. The holomorphic map
\begin{equation}
\rho(z):=(z_1+\epsilon \sum_{j,k=2}^n z_jz_k f_{jk}(z),0,\ldots, 0)
\end{equation}
is a holomorphic retraction of $\B^n$ onto the complex geodesic
$\varphi(\zeta)=(\zeta,0,\ldots, 0)$. Indeed, it is clear that
$\rho(\B^n)\subset \C\times\{(0,\ldots, 0)\}$, that $\rho^2=\rho$
and that $\rho$ is holomorphic. Moreover, if we let $r=|z_1|$ then
$|z_j|\leq \sqrt{1-r^2}$ and $|z_1+\epsilon \sum_{j,k=2}^n z_jz_k
f_{jk}(z)|\leq r+n\epsilon(1-r^2)$, proving that for $\epsilon<1/2n$
the image $\rho(\B^n)\subset \B^n$.
\end{example}

From \eqref{ro-duale} it follows that the fibers of Lempert's
projection are intersections of $D$ with complex affine hyperplanes.
Lempert's projection can be characterized exactly by this property:

\begin{prop}
Let $\varphi:\D\to D$ be a complex geodesic. If $\rho:D\to
\varphi(\D)$ is a holomorphic retraction whose fibers are
intersections of $D$ with complex affine hyperplanes then $\rho$ is
the Lempert projection. In other words, the Lempert projection is
the only ``linear'' retraction.
\end{prop}
\begin{proof}
Let $\rho:D\to \varphi(\D)$ be a retraction whose fibers are
 intersection of $D$ with complex affine hyperplanes. Let
$E_D=E_D(\varphi(e^{it}),\varphi(\zeta_0), R)$ be a horosphere of
$D$ with radius $R>0$. Since $\rho \circ \varphi=\varphi$, if $z\in
E_D$ we have
\begin{equation}\label{ins2}
\begin{split}
\lim_{w\to \varphi(e^{it})}[k_D(\rho(z),
w)&-k_D(\varphi(\zeta_0),w)]\\&=\lim_{r\to
1}[k_D(\rho(z), \rho(\varphi(re^{it})))-k_D(\varphi(0),\varphi(re^{it}))]\\
&\leq \lim_{r\to 1}[k_D(z,
\varphi(re^{it}))-k_D(\varphi(0),\varphi(re^{it}))]<\frac{1}{2}\log
R.
\end{split}
\end{equation}
Therefore $\rho(z)\in E_D\cap \varphi(\D)$. Let $\eta\in \D$ and
$\varphi(\eta)\in \de E_D$. Let $H$ be the affine hyperplane which
contains $\rho^{-1}(\varphi(\eta))$. Then $E_D\cap H=\emptyset$,
because if $z\in E_D\cap H$ then $\varphi(\eta)=\rho(z)\in E_D\cap
\varphi(\D)$, which is a contradiction. Since
$\varphi(\eta)\in\overline{E_D}\cap H$ and $E_D$ is convex, it
follows that $H-\varphi(\eta)=T^{\C}_{\varphi(\eta)}\de E_D$. Now,
$T^{\C}_{\varphi(\eta)}\de E_D=\ker (\de \rho_L)_{\varphi(\eta)}$,
where $\rho_L$ is the Lempert projection. Thus $\rho$ and $\rho_L$
have the same fibers at $\varphi(\eta)$, and, by the arbitrariness
of the choices it follows that $\rho=\rho_L$ as claimed.
\end{proof}

Next we examine the variation of the left inverse of Lempert's
projection with respect to boundary data.

\begin{lemma}\label{no-tg}
Let $\{z_k\}\subset D$ be a sequence converging non-tangentially to
$p$. Let $v_k\in L_p$ be such that $z_k\in\v_{v_k}(\D)$ (where, for
$v\in L_p$, $\v_v:\D\to D$ denotes the unique complex geodesic in
the Chang-Hu-Lee normal parametrization with respect to $v$). If
$v_{t_k}\to v_0$ then $v_0\in L_p$ and $\v_{v_k}\to \v_{v_0}$,
$\v_{v_k}^{(j)}\to \v_{v_0}^{(j)}$ uniformly on $\oD$ for all
$j=1,2,\ldots$ .
\end{lemma}
\begin{proof}
We can assume that $\nu_p=e_1$. To see that $v_0\in L_p$ we need to
show that $\la v_0, e_1\ra
>0$. Assume this is not the case. Then $v_0\in T^{\mathbb C}_p \de \D$.

First of all, we claim that  for any open neighborhood $U$ of $p$ it
follows that $\v_{v_k}(\overline{\D})\subset U$ eventually. Indeed,
let $\Phi_p:D\to \B^n$ be the spherical representation of
Chang-Hu-Lee and denote by $\eta_{v_k}:=\Phi_p \circ \v_{v_k}$. By
construction $\eta_{v_k}(\zeta)=e_1+(\zeta-1)\la v_k, e_1\ra v_k$
and thus $\eta_{v_k}(\oD)\to e_1$. Since $\Phi_p^{-1}$ is uniformly
continuous on $\overline{\B^n}$ the claim follows.

Therefore, $\{\v_{v_k}(\oD)\}$ converges to $\{p\}$ and, by
\cite[Theorem 2]{Hu}, given any $\epsilon>0$ there exists $k_0$ such
that, for all $k>k_0$, it follows that $\|(\v_{v_k}'(\zeta))_N\|\leq
\epsilon\|(\v_{v_k}'(\zeta))_T\|$ for all $\zeta\in \D$ where, if
$z\in D$ and $z'\in \de D$ is the unique point of $\de D$ nearest to
$z$, then, for all vectors $w\in T_pD=\C^n$ the vectors $w_N$ and
$w_T$ denote the complex normal and the complex tangential
components of $w$ at $z'$ (namely, $w_T\in T^{\mathbb C}_{z'}\de D$
and $w_N=\la w, \overline{\nu_{z'}}\ra \nu_{z'}$ with $\nu_{z'}$
being the unit outward normal to $\de D$ at $z'$).

Let $K\subset D$ be a cone with vertex $p$ such that $\{z_k\}\subset
K$. In particular, there exists $c>0$ such that if $w\in \C^n$ and
$(w-p)\in K$ then $\|w-p\|_N\geq c \|w-p\|_T$ (at $p$). Therefore,
if $\gamma:[0,1]\to D\cup \{p\}$ is a $C^\infty$ curve such that
$\gamma'(1)=p$ and $\|\gamma'(1)_N\|\leq (c/2)\|\gamma'(1)_T\|$ (at
$p$), then $\gamma(t)\not\in K$ for $t\approx 1$. Moreover, we can
find a small open neighborhood $U$ of $p$ such that, if
$\gamma([0,1))\subset U\cap D$ and  $\|\gamma'(t)_N\|\leq
(c/2)\|\gamma'(t)_T\|$ for all $t\in [0,1]$ (here the projection is
at the  point of $\de D$ nearest to $\gamma(t)$) then
$\gamma(t)\not\in K$ for $t\in [0,1)$).

Now, let $k$ be such that $\v_{v_k}(\D)\subset U\cap D$ and
$\|(\v_{v_k}'(\zeta))_N\|\leq (c/2)\|(\v_{v_k}'(\zeta))_T\|$ for all
$\zeta\in \D$. Let $\theta_k\in \Aut$ be an automorphism such that
$\theta_k(1)=1$ and $\v_{v_k}(\theta_k(0))=z_k$. By the previous
argument $\gamma(t):=\v_{v_k}(\theta_k(t))$ does not belong to $K$
for any $t$, which contradicts the fact that $z_k\in K$. Thus $\la
v_0, e_1\ra >0$ and $v_0\in L_p$.

We are left to show that $\v_{v_k}\to \v_{v_0}$ and
$\v_{v_k}^{(j)}\to \v_{v_0}^{(j)}$ uniformly on $\oD$. Let
$\eta_{v_k}:=\Phi_p \circ \v_{v_k}:\D\to \B^n$. By the very
definition $\eta_{v_k}(\zeta)=e_1+(\zeta-1)\la v_k,e_1\ra v_k$ and
clearly $\eta_{v_k}\to \eta_{v_0}$ uniformly on $\oD$. Since
$\Phi_p$ is a homeomorphism between $\overline{D}$ and
$\overline{\B^n}$ it follows  that $\v_{v_k}\to \v_{v_0}$ uniformly
on $\oD$. By Corollary~\ref{appendeprop} then $\v_{v_k}^{(j)}\to
\v_{v_0}^{(j)}$ uniformly on  $\oD$.
\end{proof}

\begin{lemma}\label{continuous}
For any $v\in L_p$ denote by $\v_v:\D\to D$ the unique complex
geodesic in the Chang-Hu-Lee normal parametrization with respect to
$v$ and let $\tr_v$ be its left-inverse. Then, if $\{v_k\}\subset
L_p$ is such that $v_k\to v_0\in L_p$ it follows that $d\tr_{v_k}\to
d\tr_{v_0}$ uniformly on $\overline{D}$.
\end{lemma}

\begin{proof}
Differentiating \eqref{ro-duale} with respect to $z_j$  we obtain
for $z\in \overline{D}$
\[
\frac{\de \tr_v}{\de z_j}\tv_v'(\tr_v(z))\cdot
(z-\v_v(\tr_v(z)))+\tv_v(\tr_v(z))\cdot(e_j-\frac{\de \tr_v}{\de
z_j}\v_v'(\tr_v(z)))\equiv 0,
\]
holding for $z\in \overline{D}$. Taking into account that
$\tv(\zeta)\cdot \v'(\zeta)\equiv 1$, we have
\begin{equation}\label{deride}
\frac{\de \tr_v}{\de
z_j}[\tv_v'(\tr_v(z))\cdot(z-\v_v(\tr_v(z)))-1]\equiv
-\tv_v(\tr_v(z))\cdot e_j.
\end{equation}
Notice that, since $\tv_v(\zeta)\neq 0$ for all $\zeta\in \oD$, for
all $z\in \overline{D}$ there exists $j\in \{1,\ldots, n\}$ such
that $\tv_v(\tr_v(z))\cdot e_j\neq 0$. In particular it follows that
\[
\tv_v'(\tr_v(z))\cdot(z-\v_v(\tr_v(z)))-1\neq 0
\]
for all $z\in \overline{D}$. Therefore
\begin{equation}\label{der-rho-ti}
\frac{\de \tr_v}{\de z_j}(z)= \frac{-\tv_v(\tr_v(z))\cdot e_j}{
\tv_v'(\tr_v(z))\cdot(z-\v_v(\tr_v(z)))-1}.
\end{equation}
Let $\{v_k\}\subset L_p$ be such that $v_k\to v_0\in L_p$. We claim
that
\[
\tr_{v_k}\to \tr_{v_0},\quad  \tv_{v_k}\to\tv_{v_0}, \quad
\tv_{v_k}'\to \tv_{v_0}' \quad \v_{v_k}\to \v_{v_0}
\]
uniformly on $\overline{D}$ and $\oD$ respectively. By Lemma
\ref{no-tg} it follows that $\v_{v_k}\to \v_{v_0}$  uniformly on
$\oD$.

As for $\tv_v$, if $v_k\to v_0$ in $L_p$ then by Lemma \ref{no-tg}
it follows that $\v^{(j)}_{v_k}\to \v^{(j)}_{v_0}$ uniformly on
$\oD$ for all $j=0,1,2,\ldots$. By the very definition and by
\eqref{norm-1}, if $r$ is a defining function for $\de D$, it
follows that for $\zeta\in \de\D$
\begin{equation}\label{sulbordo}
\tv_{v_k}(\zeta)=\frac{1}{\de
r_{\v_{v_k}(\zeta)}(\v_{v_k}'(\zeta))}\de r_{\v_{v_k}(\zeta)}
\end{equation}
and therefore, since $|\de
r_{\v_{v_k}(\zeta)}(\v_{v_k}'(\zeta))|\geq c>0$ for all $k$, it
follows that $\tv_{v_k}\to \tv_{v_0}$ uniformly on $\de \D$. By the
maximum principle then $\tv_{v_k}\to \tv_{v_0}$ uniformly on $\oD$.
Differentiating \eqref{sulbordo} for $\zeta=e^{it}$ and $t\in \R$ by
$\frac{d}{dt}$ we see that $\tv_{v_k}'$ is expressed as continuous
combination of $\v_{v_k}, \v_{v_k}', \v_{v_k}''$ and by Lemma
\ref{no-tg} it follows then that $\tv_{v_k}'\to \tv_{v_0}'$
uniformly on $\oD$.

We are left to show that $\tr_{v_k}\to \tr_{v_0}$ uniformly on
$\overline{D}$. If not, there exists a sequence $\{z_{k_m}\}\subset
\overline{D}$ (which we may assume converging to some $z_0\in
\overline{D}$) such that
$|\tr_{v_{k_m}}(z_{k_m})-\tr_{v_{0}}(z_{k_m})|>\epsilon_0$ for some
$\epsilon_0>0$ and for all $k_m$. By \eqref{ro-duale} it follows
that for all $k_m$
\[
\tv_{v_{k_m}}(\tr_{v_{k_m}}(z_{k_m}))\cdot
(z_{k_m}-\v_{v_{k_m}}(\tr_{v_{k_m}}(z_{k_m})))=0.
\]
Up to subsequences, we can assume that $\tr_{v_{k_m}}(z_{k_m})\to
\zeta_0\in \oD$. For what we already proved it follows then
\[
\tv_{v_0}(\zeta_0)\cdot (z_0-\v_{v_0}(\zeta_0))=0.
\]
This implies that $\zeta_0=\tr_{v_0}(z_0)$, since the only zero of
the function $\zeta\mapsto \tv_{v_0}(\zeta)\cdot
(z_0-\v_{v_0}(\zeta))$ is $\tr_{v_0}(z_0)$ by Remark \ref{zero-ro}.
But then both $\{\tr_{v_{k_m}}(z_{k_m})\}$ and
$\{\tr_{v_{0}}(z_{k_m})\}$ converge to $\tr_{v_0}(z_0)$, and then
$|\tr_{v_{k_m}}(z_{k_m})-\tr_{v_{0}}(z_{k_m})|\to 0$, contradiction.
Thus $\tr_{v_k}\to \tr_{v_0}$ uniformly on $\overline{D}$ and the
claim is proved.

Since, as we remarked at the beginning, the denominator of the right
hand side of \eqref{der-rho-ti} for $v=v_0$ is never zero for all
$z\in \overline{D}$,  the previous claim implies that $d\tr_{v_k}\to
d\tr_{v_0}$ uniformly on $\overline{D}$.
\end{proof}

\begin{nota}
By \eqref{der-rho-ti} it follows that $d(\tr_v)_p=\tv(1)$ and by
\eqref{sulbordo} we have (cfr. \cite[Lemma 2.6.44]{Ab}) for $w\in
\C^n$
\begin{equation}\label{marcolino}
d(\tr_v)_p (w)=\frac{\de r_p (w)}{\de r_p(\v_v'(1))}=\frac{\la w,
\nu_p\ra}{\la \v_v'(1), \nu_p\ra}.
\end{equation}
\end{nota}

\section{The shape of horospheres}

Let $D\subset \C^n$ be a bounded strongly convex domain and let
$p\in \de D$. As we recalled in Section \ref{prelimino}, for any
$R>0$ and $z_0\in D$, the set $\de E_D(p,z_0,R)$ is smooth away from
its center $p\in \de D$. It should be noted that smoothness of
horospheres away from the center was known after \cite[Section
4]{Trap}, but we do not know any previous reference for this fact.

In \cite{Ab} (see also \cite{Ab3}) it is proved that horospheres are
convex domains (since they are increasing union of Kobayashi balls
of $D$). In \cite[Remark 4.2]{B-P}, referring to \cite[Corollary
2.6.49]{Ab} it was claimed that (boundaries of) horospheres are
strongly convex at their center. Unfortunately the proof of
\cite[Corollary 2.6.49]{Ab} does not seem to show smoothness at the
center and thus one can only infer that horospheres are
geometrically strictly convex ({\sl i.e.}, the intersection of their
closure with the supporting hyperplane at the center is just the
center). However from \cite[p. 231-232]{Ab2} it follows that if
$E_D(p,z_0,R)\subset D$ is a horosphere of center $p\in \de D$ and
radius $R>0$ and $\B\subset D$ is a ball tangent to $\de D$ at $p$
then there exists a horosphere $E_\B(p, R')\subset \B$ for some
$R'>0$ such that $E_\B(p,R')\subset E_D(p,z_0,R)$. In particular,
since horospheres of the ball $\B$ are smooth complex ellipsoids, it
follows that there exists a ball $\B'\subset E_D(p,z_0,R)$ tangent
to $\de E_D(p,z_0,R)$ at $p$. Namely, horospheres have the
inner-ball property at the center. Therefore $\de E_D(p,z_0,R)$ is
$C^{1,1}$ at $p$ (see, {\sl e.g.}, \cite[Proposition 2.4.3]{Ho}).

We prove here that the boundaries of horospheres are strongly convex
away from the center:

\begin{teo}\label{strettaconv}
Let $D\subset \C^n$ be a strongly convex domain with smooth
boundary. Let $p\in \de D$. Let $E_D(p,R)$ be a horosphere in $D$
with center $p$ and radius $R>0$. The boundary $\de
E_D(p,R)\setminus\{p\}$ is smooth and strongly convex.
\end{teo}
\begin{proof}
Let $\Omega_{D,p}$ be the function defined in Theorem
\ref{filo-giorgio}. Its level sets are boundaries of horospheres of
$D$ with center $p$. Thus, to show that such boundaries are strongly
convex we need to prove that the (real) Hessian of $\Omega_{D,p}$ is
positive definite on the tangent space of $\de E_D(p,R)$ (for all
$R>0$). It is known (see, \cite{Ab}) that $\de E_D(p,R)$ are convex
for all $R>0$, (and  strongly pseudoconvex for all $R>0$ and
strongly convex for big radii, see \cite[Remark 7]{B-P}). Thus the
real Hessian of $\Omega_{D,p}$ is non-negative definite on the
(real) tangent space of $\de E_D(p,R)$ for all $R>0$.

Let $q\in D$ and let $\v:\D\to D$ be a complex geodesic  such that
$\v(0)=q$ and $\v(1)=p$. Up to post-composing with automorphisms of
$\B^n$ fixing $e_1$, we can suppose that $\Phi_p(q)=O$. Thus
$\Phi_p(\v(\zeta))=(\zeta,O)$. Let $F: D \to \Hr^n:=\{(\zeta,w)\in
\C\times \C^{n-1}: \Im \zeta>\|w\|^2\}$ be given by $F=C\circ
\Phi_p$ where $C:\B^n\to \Hr^n$ is the Cayley transform defined as
$$
C(\zeta,w)=(i\frac{1+\zeta}{1-\zeta}, \frac{w}{1-\zeta}),\quad
(\zeta,w)\in\C\times\C^{n-1}.
$$
We write $F(z)=(F_0(z), \tF(z))\in \C\times \C^{n-1}$. By
definition,
\begin{equation}\label{costru}
F_0(\v(\zeta))=i\frac{1+\zeta}{1-\zeta},\quad \tF(\v(\zeta))\equiv
O.
\end{equation}
By the very definition of $\Omega_{D,p}$ (see Remark \ref{comefo})
it follows that $\Omega_{D,p}(F^{-1}(\zeta,w))=\|w\|^2-\Im \zeta$
for $(\zeta,w)\in \C\times \C^n$, $(\zeta,w)\in \Hr^n$. Therefore
\begin{equation}\label{ubp}
\Omega_{D,p}(z)=\Omega_{D,p}(F^{-1}(F(z)))=\|\tF(z)\|^2-\Im F_0(z).
\end{equation}

Thus, from \eqref{costru} and \eqref{ubp} we have for $v\in \C^n$
\begin{equation}\label{hessiano}
\Hess(\Omega_{D,p})_{\v(\zeta)} (v,v)=2\|d
\tF_{\v(\zeta)}(v)\|^2-\Hess(\Im F_0)_{\v(\zeta)}(v,v)
\end{equation}
where, for a real function $f$, $\Hess (f)_x$ denotes the real
Hessian of $f$ at $x$.

Now, let $r\in \D\cap \R$ and let $\theta_r\in \Aut$ be such that
$\theta_r(0)=r$ and $\theta_r(1)=1$ (notice that necessarily
$\theta_r'(0)\in \R$). Let $\widetilde{\v\circ\theta_r}:\D\to D$ be
the dual map of $\v\circ \theta_r$. from the very definition, a
direct computation  shows that
\[
\widetilde{\v\circ\theta_r}(\zeta)=\frac{\tv(\theta_r(\zeta))}{\theta_r'(\zeta)}.
\]
By \cite[Lemma 4.1]{Trap} (and since $\theta_r'(0)\in \R$) it
follows  that if $\v(r)\in \de E_D(p,R(r))$ then
\begin{equation}\label{tangenteoro}
\begin{split}
T_{\v(r)}(\de E_D(p,R(r)))&=\{v\in \C^n : \Re (
\widetilde{\v\circ\theta_r}(0)\cdot v)=0\}\\&=\{v\in \C^n : \Re (
\widetilde{\v}(r)\cdot v)=0\}.
\end{split}
\end{equation}
On the other hand, by \eqref{ubp} and \eqref{costru} it follows that
\[
T_{\v(r)}(\de E_D(p,R(r)))=\ker d(\Omega_{D,p})_{\v(r)}=\ker d(\Im
F_0)_{\v(r)}.
\]
Thus, since they have the same kernel, the two (real) forms
$v\mapsto \Re ( \widetilde{\varphi}(r)\cdot v)$ and $v\mapsto \Im d(
F_0)_{\v(r)}(v)$ are multiple each other. Since $\Re
(\widetilde{\v}(r)\cdot \v'(r))=1$ by \eqref{norm-1}, and
by~\eqref{costru}
\begin{equation*}
\Im d(F_0)_{\varphi(r)}(\v'(r))=\Im \frac{d }{d
\xi}F_0(\v(\xi))|_{\xi=r}=\Re \frac{d }{d \xi}
\frac{1+\xi}{1-\xi}|_{\xi=r}=  \frac{2}{(1-r)^{2}},
\end{equation*}
it follows that for all $v\in\C^n$
\begin{equation}\label{d-im}
d(\Im F_0)_{\v(r)}(v)=\frac{2  \Re (\tv(r)\cdot v)}{ (1-r)^{2}}.
\end{equation}

Now, let $R>0$ be such that  $q\in \de E_D(p,R)$ and assume that
$v\in T_{q}\de E_D(p,R)$  verifies $\Hess (\Omega_{D,p})_q(v,v)=0$.
We want to show that  $v=0$.

Write $(\lambda, U)=(d (F_0)_q (v), d (\tF)_q (v))$. Since  the map
$\Phi_p$ transforms boundaries of horospheres onto boundaries of
horospheres, it follows that the vector $(\lambda, U)$ is tangent to
the boundary of the horosphere $\{(\zeta,w)\in \C\times\C^{n-1}: \Im
\zeta -\|w\|^2>1\}$ whose closure contains $(i,O)\in \Hr^n$. Thus
$\lambda\in \R$.

Let us now consider the smooth one-parameter family of complex
geodesics  $g_t:\D\to \Hr^n$ depending smoothly on $t$ given by
\[
g_t(\zeta):=(i\frac{1+\zeta}{1-\zeta}+t\lambda+it^2\|U\|^2, tU)
\]
and we denote
\[
\gdot (\zeta):=\frac{\de g_t(\zeta)}{\de t}|_{t=0}, \quad \gddot
(\zeta):=\frac{\de^2 g_t(\zeta)}{\de^2 t}|_{t=0}.
\]
Notice that $\gdot (\zeta)=(\lambda, U)$ and
$\gddot(\zeta)=(2i\|U\|^2,O)$  are independent of $\zeta$. Let
$f_t:=F^{-1}(g_t)$. By construction $f_t:\D\to D$ is a smooth
one-parameter family of complex geodesics, $f_t(1)=p$ and
$f_0(0)=q$. Therefore $f_0=\v$. Again,  denoting by $\fdot(\zeta),
\fddot(\zeta)$ the derivative of $f_t(\zeta)$ with respect to $t$ at
$t=0$, it follows that $\fdot(1)=0$ and $\fddot(1)=0$ because
$f_t(1)=p$ for all $t$ (see Corollary \ref{appende}).

Let us denote by $J(\zeta)$ the {\sl Jacobi vector field
$J(\zeta):=\fdot(\zeta)$ along $\v$}. We can write
$J(\zeta)=\lambda(\zeta)\v'(\zeta)+J^\perp(\zeta)$ for some
holomorphic function $\lambda$ and vector field $J^\perp$ such that
$\tv(\zeta)\cdot J^\perp(\zeta)\equiv 0$ and
$\lambda(\zeta)=\alpha+i\beta\zeta-\overline{\alpha}\zeta^2$ for
some $\alpha\in \C$ and $\beta\in \R$ (see \cite[Theorem
3.1.c]{Trap}).

Since $\fdot(1)=0$ then $J(1)=0$. Since the map $\Phi_p$ transforms
boundaries of horospheres onto boundaries of horospheres, it follows
that $J(0)\in T_{\v(0)}(\de E_D(p,R))$. In other words, by
\eqref{tangenteoro}, $\Re (\tv(0)\cdot J(0))=0$ which implies that
$\Re\alpha=0$ for $\tv(0)\cdot \v'(0)=1$. Therefore
$J(\zeta)=i\gamma (1-\zeta)^2\v'(\zeta)+J^\perp(\zeta)$ with
$\gamma\in \R$ and $J^\perp(1)=0$. This implies that $\Re
(\tv(r)\cdot J(r))=0$ for $r\in (-1,1)$. Hence by
\eqref{tangenteoro} it follows that $J(r)\in T_{\v(r)}(\de
E_D(p,R(r)))$, where $R(r)>0$ is such that $\v(r)\in \de E_D(p,
R(r))$. Since boundaries of horospheres are convex, we have
\begin{equation}\label{Jmag}
\Hess(\Omega_{D,p})_{\v(r)}(J(r), J(r))\geq 0,\quad r\in (-1,1).
\end{equation}

Now differentiating with respect to $t$ the identity $F \circ f_t=
g_t$ and setting $t=0$ we obtain
\begin{equation}\label{diff-rent}
\begin{split}
&d \tilde{F}_{\v(\zeta)}(J(\zeta))= U \\
& \Hess (\Im F_0 )_{\v(\zeta)}(J(\zeta), J(\zeta))+ d (\Im
F_0)_{\v(\zeta)}(\fddot(\zeta))=\Im \gddot (\zeta)=2\|U\|^2.
\end{split}
\end{equation}
Putting together \eqref{hessiano}, \eqref{d-im} and
\eqref{diff-rent}, we obtain for $r\in (-1,1)$
\begin{equation}\label{ess-lungor}
\Hess (\Omega_{D,p})_{\v(r)}(J(r),J(r))=\frac{\Re\tv(r)\cdot\fddot
(r) }{(1-r)^2}.
\end{equation}

Our next aim is to compute $\tv(r)\cdot\fddot (r)$. In order to do
this, we choose a suitable defining function: according to Pang
\cite[Proposition 2.36]{Pg1} there exists a $C^\infty$ defining
function $\rho$ for $D$ near $\v(\overline{D})$ such that for all
$\theta\in \R$ it follows that $\tv(e^{i\theta})=e^{i\theta}\de
\rho_{\v(e^{i\theta})}$. For all $t$ and for all $\theta\in \R$ it
follows that $\rho(f_t(e^{i\theta}))\equiv 0$, thus differentiating
such an identity (as we can, by Corollary \ref{appende}) with
respect to $t$ at $t=0$ we obtain $2\Re (\de \rho \cdot
\fddot(e^{i\theta})+\Hess (\rho)_{e^{i\theta}}(J(e^{i\theta}),
J(e^{i\theta}))\equiv 0$, namely,
\begin{equation}\label{b1}
\Re (\overline{\zeta} \tv(\zeta)\cdot \fddot
(\zeta))=-\frac{1}{2}\Hess (\rho)_{\zeta}(J(\zeta), J(\zeta)), \quad
|\zeta|=1.
\end{equation}
Now, the function $\zeta\mapsto \tv(\zeta)\cdot \fddot(\zeta)$ is
holomorphic in $\D$. Write $\tv(\zeta)\cdot \fddot(\zeta)=A+\zeta
B+\zeta^2 C(\zeta)$ for some $A, B\in \C$ and some holomorphic
function $C$. Then
\begin{equation}\label{b2}
\Re (\overline{\zeta} \tv(\zeta)\cdot \fddot (\zeta))=\Re
(\overline{A}\zeta+B+\zeta C(\zeta)),\quad |\zeta|=1.
\end{equation}
Let $T_1$ denote the Hilbert transform which associates to any
harmonic function $u$ in $\D$, H\"older continuous on $\de \D$, its
harmonic conjugated $T_1(u)$, still H\"older continuous on $\de \D$,
normalized so that $T_1(u)(1)=0$. Let $h$ denote the holomorphic
function in $\D$ whose trace on $\de \D$ is $-1/2({\sf id} + i
T_1)(\Hess (\rho)(J,J))$. Notice that $\Re h\leq 0$ on $\de \D$
since $\de D$ is convex. Moreover, since $J(1)=0$ and by the
normalization chosen for $T_1$ it follows that $h(1)=0$.

By \eqref{b1} and \eqref{b2} we obtain that
$h(\zeta)=\overline{A}\zeta+B+\zeta C(\zeta)+i\alpha$ for some
$\alpha\in \R$. Hence $\tv(\zeta)\cdot \fddot(\zeta)-\zeta
h(\zeta)=-\overline{A}\zeta^2-i\alpha \zeta +A$ and, since
$h(1)=\fddot(1)=0$, we get $-\overline{A}-i\alpha  +A=0$. Writing
$A=a+ib$ with $a,b\in \R$, we obtain
\[
\tv(\zeta)\cdot \fddot(\zeta)=\zeta
h(\zeta)+a(1-\zeta^2)+ib(1-\zeta)^2.
\]
Substituting this expression in \eqref{ess-lungor} for $\zeta=r\in
(-1,1)$, we find
\begin{equation}\label{finalex}
\Hess (\Omega_{D,p})_{\v(r)}(J(r),J(r))=a\frac{1+r}{1-r}+\frac{r\Re
(h(r))}{(1-r)^2}.
\end{equation}
By construction then
$a=\Hess(\Omega_{D,p})_{\v(0)}(J(0),J(0))=\Hess(\Omega_{D,p})_{q}(v,v)=0$.
By \eqref{Jmag} and \eqref{finalex} it follows then that $\Re
(h(r))\geq 0$ for $r\in (0,1)$. However $\Re h(\zeta)$ is harmonic
on $\D$ and non-positive on $\de \D$ and thus by the maximum
principle $\Re h(\zeta)\equiv 0$. Thus
\[
\Hess (\rho)_\zeta (J(\zeta), J(\zeta))= 0, \quad \zeta \in \de \D
\]
and, since $\de D$ is strongly convex, it follows that $J=0$ on $\de
\D$ and thus $J\equiv 0$ on $\D$ proving that $v=0$.
\end{proof}

\section{Extremality}\label{estremalita}

Let $D\subset \C^n$ be a bounded strongly convex domain with smooth
boundary. We let $\Gamma_p$ be the set of all $C^\infty$ curves
$\gamma:[0,1]\to D\cup\{p\}$ such that $\gamma(1)=p$ and
$\gamma'(1)\not\in T_p\de D$ (notice that, if $\nu_p$ is the unit
outward normal to $\de D$ at $p$ then $\gamma'(1)\not\in T_p\de D$
if and only if $\Re\la \gamma'(1),\nu_p\ra> 0$).

\begin{teo}\label{extremol}
Let $D\subset \C^n$ be a bounded strongly convex domain with smooth
boundary and let $p\in \de D$. Let $\nu_p$ be the unit outward
normal to $\de D$ at $p$. Consider the following family $\mathcal
S_{p}(D)$:
\begin{equation}\label{family1}
\begin{cases}
u \in \ps (D) \\
\limsup_{z\to x} u(z)\leq  0 \quad \hbox{for all}\ x\in \de D\setminus\{p\}\\
\displaystyle{\liminf_{t\to 1} |u(\gamma(t))(1-t)|\geq 2\Re(\la
\gamma'(1),\nu_p\ra^{-1})} \quad \hbox{for all }\gamma\in \Gamma_p,
\end{cases}
\end{equation}
Then $\Omega_{D,p}\in \mathcal S_{p}(D)$ (where $\Omega_{D,p}$ is
the function defined in Theorem \ref{filo-giorgio}) and $u\leq
\Omega_{D,p}$ for all $u\in \mathcal S_{p}(D)$.
\end{teo}

To prove the theorem we need some preliminary results.  Let
$\subh(\D)$ denote the real cone of subharmonic functions in the
unit disc $\D$.

\begin{lemma}[Phragmen-Lindel\"of]\label{uno}
Let $c>0$. Consider the following family $\mathcal S_{c}(\D)$ in the
unit disc:
\begin{equation}\label{family-disc}
\begin{cases}
u \in \subh (\D) \\
u<0 \quad \hbox{in}\ \D \\
\displaystyle{\liminf_{\R\ni r\to 1^-} |u(r)(1-r)|\geq 2c}
\end{cases}
\end{equation}
Then $-cP(\zeta)\in \mathcal S_c(\D)$ and for all $u\in \mathcal
S_c(\D)$ it follows $u\leq -cP(\zeta)$.
\end{lemma}

For the sake of completeness we give a short proof of Lemma
\ref{uno}, based on some notes of Prof. P. Poggi-Corradini. We thank
him for letting us to use such notes.
\begin{proof}
It is clear that $-cP(\zeta)\in \mathcal S_c(\D)$. We have to show
that $-cP$ is the maximal element of the family.

First of all, let $C(\zeta)=(1+\zeta)\cdot(1-\zeta)^{-1}$ be the
Cayley transformation from $\D$ to $\Hr=\{w\in \C: \Re w>0\}$. Then
we consider the family $C^\ast(\mathcal S_c(\D)))=\{\tilde{u} :
\tilde{u}=u\circ C^{-1} \hbox{ for some } u\in \mathcal S_c(\D)\}$.
Then, if $\tu\in C^\ast(\mathcal S_c(\D))$ it follows that
$\tu\in\subh(\Hr)$, $\tu<0$ in $\Hr$ and
\[
\limsup_{\R\ni R\to +\infty} \frac{\tu(R)}{R}\leq -c.
\]
Notice that $P\circ C^{-1}(w)=\Re w$ is the  Poisson kernel in
$\Hr$. Let $\tilde{u}=u\circ C^{-1}\in C^\ast(\mathcal S_c(\D))$ and
let $L= \limsup_{\R\ni R\to +\infty} \tu(R)/R\leq -c$. We are going
to show that $\tu\leq L\Re w$, from which it follows that $u\leq
-cP$.

Let $\epsilon>0$ be such that $\epsilon<-L$. Let
$v(w)=\tu(w)-(L+\epsilon)\Re w$. Now, $v\in\subh(\Hr)$,
$\limsup_{w\to iy}v(w)\leq 0$ for all $y\in \R$ and
\[
\limsup_{\R\ni R\to +\infty}v(R)=\limsup_{\R\ni R\to
+\infty}R\left(\frac{\tu (R)}{R}-(L+\epsilon)\right)\leq 0.
\]
Therefore there exists $\delta>0$ such that $v(R)\leq 1$ for $R\leq
\delta$ and $R\geq \frac{1}{\delta}$. Moreover, since $v$ is
semicontinuous, there exists $K>0$ such that $v(R)\leq K$ for
$\delta<R<\frac{1}{\delta}$. We consider now $V(w)=v(\sqrt{iw})-K$.
Again, $V\in \subh(\Hr)$ and $\limsup_{w\to iy}V(w)\leq 0$ for all
$y\in \R$. Moreover,
\begin{equation*}
\begin{split}
\sup_{-\pi/2<\theta<\pi/2}
V(re^{i\theta})&=\sup_{0<\theta<\pi/2}v(\sqrt{r}e^{i\theta})-K=\sup_{0<\theta<\pi/2}[\tu(\sqrt{r}e^{i\theta})
-(L+\epsilon)\sqrt{r}\cos \theta-K]\\ &\leq
\sup_{0<\theta<\pi/2}(-(L+\epsilon)\sqrt{r}\cos
\theta-K)=-(L+\epsilon)\sqrt{r}-K.
\end{split}
\end{equation*}
By the classical estimates on sub-linear growth of subharmonic
functions (see, {\sl e.g.} \cite{P-W}), it follows that $V(w)\leq 0$
for all $w\in \Hr$ and therefore, $v\leq K$ in the first quadrant. A
similar argument shows that $v\leq K$ on the fourth quadrant and as
before, $v\leq 0$ on $\Hr$ which implies $\tu(w)\leq (L+\epsilon)\Re
w$ for $w\in \Hr$. By the arbitrariness of $\epsilon$ we have the
statement.
\end{proof}

\begin{proof}[Proof of Theorem \ref{extremol}]
Up to rigid movements, we can suppose that $\nu_p=e_1$.

First of all,  notice that the function identically $0$ does not
belong to $\mathcal S_p(D)$ because of the estimates at $p$.

We claim that if $u\in \mathcal S_p(D)$ then $u<0$ in $D$. Indeed,
 let $\v:\D\to D$ be a complex geodesic  not containing $p$ in its closure (in fact, any
attached analytic disc not containing $p$ would be enough). Then
$\tu=u\circ \v:\D\to \R$ is subharmonic and $\limsup_{\zeta\to
x}\tu(\zeta)\leq 0$ for all $x\in \de \D$. Thus by the maximum
principle for subharmonic functions, $\tu\leq 0$ in $\D$ and hence
$u\leq 0$ in $D$ as $\v$ was an arbitrary complex geodesic. Again,
the maximum principle for plurisubharmonic functions implies that
$u<0$ in $D$ or $u\equiv 0$, and the latter cannot be the case. Thus
$u<0$ in $D$ as wanted.

Now, let $v\in L_p$ and let $\v_v:\D\to D$ be a complex geodesic
parameterized as in \cite{CHL}. Let $\tr_v:D\to \D$ be its left
inverse. We show that the function $u_v:D\to \R^-$ defined by
\begin{equation}\label{uv}
u_v(z)=-\frac{P(\tr_v(z))}{v_1^2 }
\end{equation}
belongs to $\mathcal S_p(D)$. It is clear that $u_v\in \ps(D)$,
$\limsup_{z\to x}u_v(z)\leq 0$ for all $x\in \de D\setminus\{p\}$
since $\tilde{\rho}_v(\overline{D}\setminus \varphi_v(\de
\D))\subset \D$. We claim that for all smooth curves
$\gamma:[0,1]\to D\cup \{p\}$ such that $\gamma(1)=p$ and $\la
\gamma'(1), e_1\ra \neq 0$ (that is $\gamma'(1)$ is not complex
tangential to $\de D$) it follows
\begin{equation}\label{uv-va}
\lim_{t\to 1} |u_v(\gamma(t))(1-t)|= \frac{2\Re\la {\gamma}'(1),
e_1\ra }{|\la {\gamma}'(1),e_1\ra|^2}.
\end{equation}
Indeed,
\[
|u_v(\gamma(t))(1-t)|=\frac{1}{v_1^2}\frac{1-|\tr_v(\gamma(t))|^{2}}{1-t}\frac{(1-t)^2}{|1-\tr_v(\gamma(t))|^2},
\]
now
\[
\lim_{t\to
1}\frac{1-\tr_v(\gamma(t))}{1-t}=\frac{d}{dt}(\tr_v(\gamma(t)))|_{t=1}=d(\tr_v)_p(\gamma'(1))
=\frac{\la\gamma'(1),e_1\ra}{\la \v_v'(1),e_1\ra},
\]
where the last equality follows from \eqref{marcolino} and since
$\v_v'(1)=v_1v$. Moreover
\begin{equation*}
\begin{split}
\lim_{t\to
1}\frac{1-|\tr_v(\gamma(t))|^2}{1-t}&=\frac{d}{dt}(|\tr_v(\gamma(t))|^2)|_{t=1}\\&=
\tr_v(\gamma(1))\frac{d}{dt}(\overline{\tr_v(\gamma(t))})|_{t=1}
+\overline{\tr_v(\gamma(1))}\frac{d}{dt}(\tr_v(\gamma(t)))|_{t=1}\\&=2\Re
\frac{d}{dt}(\tr_v(\gamma(t)))|_{t=1}=2\Re
\frac{\la\gamma'(1),e_1\ra}{\la \v_v'(1),e_1\ra}.
\end{split}
\end{equation*}
Therefore
\begin{equation*}
\lim_{t\to 1} |u_v(\gamma(t))(1-t)|=\frac{2}{v_1^2}\Re
\frac{\la\gamma'(1),e_1\ra}{\la \v_v'(1),e_1\ra}\cdot \frac{|\la
\v_v'(1),e_1\ra|^2}{|\la\gamma'(1),e_1\ra|^2}.
\end{equation*}
Taking into account that $\la \v_v'(1), e_1\ra =v_1^2$, we have the
claim. In particular equation \eqref{uv-va} holds if $\gamma\in
\Gamma_p$, showing that $u_v$ belongs to $\mathcal S_p(D)$.

Notice that $\Omega_{D,p}(\varphi_v(\zeta))=u_v(\varphi_v(\zeta))$
for all $\zeta\in \D$. Moreover, if $u\in \mathcal S_p(\D)$ then
$\tu:\zeta\mapsto u(\varphi_v(\zeta))$ is in the family $\mathcal
S_{1/v_1^2}(\D)$ given by \eqref{family-disc} (with $c=1/v_1^2$).
Indeed, $\tu\in \subh(\D)$, $u<0$ in $\D$ and
\[
\liminf_{\R\ni r\to 1}|\tu(r)|(1-r)\geq \frac{2\Re\la
\v_v'(1),e_1\ra}{|\la\v_v'(1),e_1\ra|^2}=\frac{2}{\la\v_v'(1),e_1\ra}=\frac{2}{v_1^2},
\]
since $\v_v'(1)=v_1^2e_1+e_1^\perp$. Thus, by Lemma \ref{uno} it
follows that for all $\zeta\in\D$
\[
u(\v_v(\zeta))=\tu(\zeta)\leq
\frac{-1}{v_1^2}P(\zeta)=\frac{-1}{v_1^2}u_v(\v_v(\zeta))=\Omega_{D,p}(\varphi_v(\zeta)).
\]
Thus for all $u\in \mathcal S_p(D)$ we have $u\leq \Omega_{D,p}$. It
remains only to show that $\Omega_{D,p}\in \mathcal S_p(D)$. To this
aim, we let $\v_{v_t}:\D\to D$ be the complex geodesic in
Chang-Hu-Lee normal parametrization such that $\gamma(t)\in
\v_{v_t}(\D)$. Moreover, we denote by $v_t=\v_{v_t}'(1)\in L_p$.
Thus
\[
\Omega_{D,p}(\gamma(t))=u_{v_t}(\gamma(t))=\frac{-1}{\la v_t,
e_1\ra^2}P(\tr_{v_t}(\gamma(t))).
\]
Hence
\begin{equation}\label{omega-bordo}
|\Omega_{D,p}(\gamma(t))|(1-t)=\frac{1}{\la v_t,
e_1\ra^2}\frac{1-|\tr_{v_t}(\gamma(t))|^2}{1-t}\frac{(1-t)^2}{|1-\tr_{v_t}(\gamma(t))|^2}.
\end{equation}
Fix $v=v_t$. By the mean value theorem it follows that
\[
\frac{1-\Re\tr_{v}(\gamma(t))}{1-t}=\frac{d}{dt}\Re\tr_{v}(\gamma(t))|_{t=s}=\Re
d(\tr_{v})_{\gamma(s)}(\gamma'(s)),
\]
for some $t<s<1$, and similarly for the imaginary part and for the
modulus $|\tr_{v_t}(\gamma(t))|^2$. Notice that $s$ depends on $v$
but clearly, $s\to 1$ as $t\to 1$.

Now let $\{v_{t_k}\}$ be a converging subsequence of $\{v_t\}$. By
Lemma \ref{no-tg} if $v_{t_k}\to v_0$ then $v_0\in L_p$ (and in
particular $\la v_{t_k},e_1\ra^2\to\la v_{0},e_1\ra^2>0$).
Therefore, by Lemma \ref{continuous} we have
\[
\lim_{t\to 1} d(\tr_{v_t})_{\gamma(s)}(\gamma'(s))=
d(\tr_{v_0})_{p}(\gamma'(1)).
\]
Thus by \eqref{omega-bordo} and  \eqref{marcolino}  it follows
\[
\lim_{{t_k}\to 1}|\Omega_{D,p}(\gamma(t_k))|(1-t_k)=\frac{1}{\la
v_0, e_1\ra^2}\frac{2\Re
d(\tr_{v_0})_p(\gamma'(1))}{|d(\tr_{v_0})_p(\gamma'(1))|^2}=\frac{2\Re\la
{\gamma}'(1), e_1\ra }{|\la {\gamma}'(1),e_1\ra|^2}.
\]
Since this holds for any converging subsequence of $\{v_t\}$ then we
have that
\[
\lim_{{t}\to 1}|\Omega_{D,p}(\gamma(t))|(1-t)=\frac{2\Re\la
{\gamma}'(1), e_1\ra }{|\la {\gamma}'(1),e_1\ra|^2}.
\]
\end{proof}

\begin{corollario}\label{oicomevobene}
Let $\Omega_{D,p}$ be the function given by Theorem
\ref{filo-giorgio}. Then for all smooth curves $\gamma:[0,1]\to
D\cup\{p\}$ such that $\gamma(1)=p$ and $\gamma'(1)\not\in T^\C_p\de
D$ it follows
\[
\lim_{t\to 1}|\Omega_{D,p}(\gamma(t))|(1-t)=\Re\frac{2}{\la
{\gamma}'(1),\nu_p\ra}.
\]
\end{corollario}
\begin{proof}
If $\gamma'(1)\not\in T_p\de D$ then the claim follows from the
proof of Theorem \ref{extremol}. In case $\gamma'(1)\in T_p\de
D\setminus T_p^\C\de D$---that is $\Re \la \gamma'(1), \nu_p \ra =0$
but $\la \gamma'(1), \nu_p \ra \neq 0$---let $v\in L_p$ and let
$u_v$ be given by \eqref{uv}. By Theorem \ref{extremol} it follows
that for all $z\in D$
\[
0\leq |\Omega_{D,p}(z)|\leq |u_v(z)|.
\]
By \eqref{uv-va} it follows that $|u_v(\gamma(t))|(1-t)\to 0$ and
then $|\Omega_{D,p}(\gamma(t))|(1-t)\to 0$, proving the statement.
\end{proof}

\section{Green's  versus Poisson's pluricomplex  functions}

Let $D$ be a bounded strongly convex domain in $\C^n$ with smooth
boundary and let $z_0\in D$. Consider the  problem in
\eqref{monge-inner-intro}. In his outstanding work \cite{Le},
\cite{Le3}, Lempert proved that there exists  a  unique solution
$L_{D,z_0}$, given by $L_{D, z_0}=\log \|\Phi_{z_0}\|$, where
$\Phi_{z_0}:D\to \B^n$ is the Lempert spherical representation with
center $z_0$ introduced in Section \ref{prelimino}. Rephrasing the
very definition of $\Phi_{z_0}$, it follows that for $z\in D$
\begin{equation}\label{lempert}
L_{D, z_0}(z)=\log (\tanh k_D(z_0, z)).
\end{equation}

 We have the following relations between the pluricomplex Green function $L_{D, z_0}$
and the pluricomplex Poisson kernel $\Omega_{D,p}$ solution of the
problem \eqref{monge-boundary} which generalizes the corresponding
relation in $\D$ between the classical Green function and the
classical Poisson kernel (see for instance \cite[Proposition
2.2.2]{Kl}):
\begin{teo}\label{GvP}
Let $D$ be a bounded strongly convex domain in $\C^n$ with smooth
boundary. Let $z_0\in D$ and $p\in \de D$. Let $\nu_p$ be the outer
normal of $\de D$ at $p$. Then
\begin{equation}\label{legati}
\Omega_{D,p}(z_0)=-\frac{\de L_{D,z_0}}{\de \nu_p}(p)
\end{equation}
\end{teo}
\begin{proof} Let $K_{z_0}:=\exp(L_{D,z_0})$.
 Let $\v:\D\to D$ be a
complex geodesic such that $\v(0)=z_0$ and $\v(1)=p$. Since
$\frac{\de K_{z}}{\de \nu_p}(p)>0$ for all $z\in D$, by
\cite[Theorem 2.6.47]{Ab} (see also \cite{Ab3}) it follows that
\[
\lim_{\R\ni t\to 1}[k_D(z,\v(t))-k_D(z_0, \v(t))]=\frac{1}{2}[\log
\frac{\de K_{z_0}}{\de \nu_p}(p)-\log \frac{\de K_{z}}{\de
\nu_p}(p)].
\]
On the other hand by \cite[Proposition 7.1]{B-P}
\[
\lim_{\R\ni t\to 1}[k_D(z,\v(t))-k_D(z_0, \v(t))]=\frac{1}{2}[\log
|\Omega_{D,p}(z_0)|-\log|\Omega_{D,p}(z)|],
\]
which implies that there exists $C>0$ such that for all $z\in D$
\[
\Omega_{D,p}(z)=-C \frac{\de K_{z}}{\de \nu_p}(p).
\]
We want to show that $C=1$. Let $\v:\D\to D$ be the unique complex
geodesic in Chang-Hu-Lee normal parametrization such that $\v(1)=p$
and $\v'(1)=\nu_p$. By the very definition
$\Omega_{D,p}(\v(\zeta))=-P(\zeta)$, where $P$ is the Poisson kernel
of $\D$ and $K_{\v(0)}(\v(\zeta))=|\zeta|$ for all $\zeta\in \D$.
Since
\[
\frac{\de K_{\v(0)}}{\de \nu_p}(p)=\frac{d}{dr}(K_{\v(0)}\circ
\v(r))|_{r=0}=\frac{d}{dr}r=1,
\]
and $P(0)=1$ it follows that $C=1$, as wanted. Finally, since
$\frac{\de K_{z}}{\de \nu_p}(p)=K_z(p)\frac{\de L_{D,z}}{\de \nu_p}$
and $K(p)=1$ for $p\in \de D$, we get \eqref{legati}.
\end{proof}

\section{Uniqueness properties}

In this section we study some  analytical and geometrical properties
which characterize the pluricomplex Poisson kernel introduced
before.

Before start, recall that, according to Bedford and Taylor
\cite{B-T2} (see also \cite[Section 3.5]{Kl}, the complex
Monge-Amp\`ere operator $(dd^c)^n$ (here $d^c=i(\debar-\de)$) can be
defined for all $u\in \ps(D)\cap L^\infty_{loc}(D)$ for any bounded
domain $D\subset \C^n$. Moreover, if $u\in \ps(D)\cap
L^\infty_{loc}(D)$ then $(dd^c u)^n=(\de \debar u)^n=0$ if and only
if $u$ is {\sl maximal} in $D$; namely, for all relatively compact
open subsets $E\subset D$ and all  plurisubharmonic functions $v$ in
$E$  such that $\limsup_{E\ni z\to x}v(z)\leq u(x)$ for all $x\in\de
E$ it follows that $v\leq u$ in $E$.

Now we can state and prove the first uniqueness result, which is the
analogous in our setting of the uniqueness statement for the
Monge-Amp\`ere equation with one concentrated logarithmic
singularity in the domain $D$ (see \cite{Le3}).

\begin{teo}\label{unico-comportamento-al-bordo}
Let $D\subset \C^n$ be a bounded strongly convex domain with smooth
boundary and let $p\in \de D$. Let $u\in\ps(D)\cap
L_{loc}^\infty(D)$ be such that $(\de \debar u)^n=0$, $\lim_{z\to
x}u(z)=0$ for all $x\in \de D\setminus\{p\}$ and
\begin{equation}\label{albordo}
\lim_{z\to p} \frac{u(z)}{\Omega_{D,p}(z)}=1.
\end{equation}
Then $u\equiv \Omega_{D,p}$.
\end{teo}
\begin{proof}
First of all we notice that \eqref{albordo} implies that $u$ belongs
to the family \eqref{family1} because for all $\gamma\in \Gamma_p$
(here $\Gamma_p$ is the set of curves defined in section
\ref{estremalita}) it follows that
\[
\lim_{t\to 1}u(\gamma(t))(1-t)=\lim_{t\to
1}\frac{u(\gamma(t))}{\Omega_{D,p}(\gamma(t))}\Omega_{D,p}(\gamma(t))(1-t)=-\Re2(\la
{\gamma}'(1),\nu_p\ra)^{-1}.
\]
Therefore, by Theorem \ref{extremol} it follows that $u(z)\leq
\Omega_{D,p}(z)$ for all $z\in D$. Suppose that
$u(z_0)<\Omega_{D,p}(z_0)$ for some $z_0\in D$. Then there exist
$0<c<1$ and $\delta>0$ such that  the set
\[
E_{\delta,c}:=\{z\in D: \Omega_{D,p}(z)>cu(z)+\delta\}
\]
is non-empty. Since $u$ is upper semi-continuous the set
$E_{\delta,c}$ is open. If we prove that $E_{\delta,c}$ is
relatively compact in $D$, since $(\de\debar (cu+\delta))^n=0$ and
$\Omega_{D,p}(z)\leq cu(z)+\delta$ on $\de E_{\delta,c}$, by
maximality it follows that $\Omega_{D,p}(z)\leq cu(z)+\delta$ in
$E_{\delta,c}$, contradicting the definition of $E_{\delta, c}$.

Thus we are left to show that $E_{\delta,c}$ is relatively compact
in $D$. First of all, since $u(x)=\Omega_{D,p}(x)=0$ for all $x\in
\de D\setminus\{p\}$, then $\overline{E_{\delta,c}}\subset
D\cup\{p\}$. Seeking  a contradiction, we assume that $p\in
\overline{E_{\delta,c}}$. Thus there exists  $\{z_k\}\subset
E_{\delta,c}$  such that $z_k\to p$. Therefore for all $k\in \N$
\begin{equation}\label{operazioneZ}
\Omega_{D,p}(z_k)-cu(z_k)-\delta>0.
\end{equation}
Up to subsequences, we can assume that $\Omega_{D,p}(z_k)\to L$ for
some $L\in [-\infty, 0]$. If $L<0$ then dividing \eqref{operazioneZ}
by $\Omega_{D,p}(z_k)<0$ and passing to the limit, taking into
account \eqref{albordo}, we would find $1-c\leq 0$, a contradiction
since $c<1$. If $L=0$, \eqref{albordo} implies that $u(z_k)\to 0$ as
$k\to \infty$ and therefore we reach a contradiction  by passing to
the limit for $k\to\infty$ in \eqref{operazioneZ}. Hence $p$ is not
in the closure of $E_{\delta,c}$ which is thus relatively compact in
$D$.
\end{proof}

\begin{nota}
As pointed out in the introduction, Theorem
\ref{unico-comportamento-al-bordo} is the analogous of the
uniqueness statement for the problem \eqref{monge-inner-intro},
where uniqueness is established in the class of plurisubharmonic
functions $u\in \ps(D)$  such that $\lim_{z\to x}u(z)=0$ for all
$x\in \de D$ and $u(z)$ goes like the pluricomplex Green function
$L_{D,z_0}$ for $z\to z_0$. Since for any convex domain the function
$L_{D,z_0}$ goes like $\log\|z-z_0\|$ at $z_0$ then in the case of a
inner singularity, there is a ``universal'' behavior. When the
singularity is at $p\in\de D$, it turns out that, thanks to
Corollary \ref{oicomevobene}, we know that the behavior of
$\Omega_{D,p}$ along all non tangential directions is independent of
$D$, but we do not have any hint on the behavior of $\Omega_{D,p}$
along the tangential directions, which might depend on $D$ near $p$.
\end{nota}

Next we characterize the pluricomplex Poisson kernel in terms of its
associated Monge-Amp\`ere foliation, with no hypotheses on the
behavior near the boundary singularity.

\begin{teo}\label{foliazione-unica}
Let $D\subset \C^n$ be a bounded strongly convex domain with smooth
boundary and let $p\in \de D$. Let $u\in\ps(D)\cap C^2(D)$ be such
$\lim_{z\to x}u(z)=0$ for all $x\in \de D\setminus\{p\}$. Then the
restriction of $u$ to each complex geodesic whose closure contains
$p$ is harmonic if and only if  there exists $c\geq 0$ such that
$u=c\Omega_{D,p}$.
\end{teo}

\begin{proof} One direction is obvious. Assume then that $u\in\ps(D)\cap C^2(D)$ is harmonic on each complex
geodesic whose closure contains $p$ and $\lim_{z\to x}u(z)=0$ for
all $x\in \de D\setminus\{p\}$. Arguing as in the proof of Theorem
\ref{extremol} we see that $u<0$ in $D$ or $u\equiv 0$. In the
latter case $c=0$ and the theorem is proved. Thus we can assume that
$u<0$ in $D$.

First of all, it is a  well known result that if $v\geq 0$ is a
harmonic function in $\D$ such that $\lim_{\zeta\to x}v(\zeta)=0$
for all $x\in \de \D \setminus\{1\}$ then $v=c P$ for some $c\geq 0$
(here, as usual, $P$ denotes the Poisson kernel).

Therefore $u=\lambda \Omega_{D,p}$ for some $C^2$ function
$\lambda:D\to (0,\infty)$ which is constant on each complex geodesic
whose closure contains $p$. We need to show that $\lambda$ is
constant.

 To this aim, we argue as in the proof of Theorem
\ref{strettaconv} and retain the notations introduced there. Let
$q\in D$. Up to post-composing with automorphisms of $\B^n$ and with
the Cayley transform, we let $F:D\to \Hr^n$ be the diffeomorphism
defined by means of the boundary spherical representation $\Phi_p$,
such that $F(q)=(i,O)$. We let $U=u\circ F^{-1}$. Then $U$ is a
$C^2$ negative function on $\Hr^n$ and by the very definition of
$\Omega_{D,p}$ and \cite[Theorem 7.3]{B-P}, it follows
$U(\xi,w)=\tl(w)(\|w\|^2-\Im \xi)$. We are going to prove that $w=O$
is a critical point for $\tl$; from this, since $F$ is a
diffeomorphism from $D$ to $\Hr^n$, it will follow that $\lambda$
has a critical point at $q=F^{-1}(i,O)$ and, by the arbitrariness of
$q$, it will follow that all points of $D$ are critical for
$\lambda$ which turns out to be  constant.

Since $\tl$ is a real function, it is enough to prove that the
vector $V:=(\frac{\de \tl}{\de w_1}(0), \ldots, \frac{\de \tl}{\de
w_{n-1}}(0))$ is zero. Let $\v:\D \to D$ be the complex geodesic
such that $\v(0)=q$ and $\v(1)=p$. According to \cite[Section
2.39]{Pg1} we can assume to be working with a system of holomorphic
coordinates $(z_1,\ldots, z_n)$ in a neighborhood  of $\v(\oD)$ for
which (among other conditions on the defining function of $D$ which
we only use implicitly when referring to the paper \cite{S-T} in the
course of the proof) $\v(\zeta)=(\zeta,0,\ldots, 0)$ for $\zeta\in
\D$.

 By construction it
follows that if we write $G:=F^{-1}=(G_1,\ldots, G_n)$ then
 $G_1(\xi,O)=(\xi-i)/(\xi+i)$
and $G_j(\xi,O)=0$ for $j>1$ and $\Im \xi>0$.

Now let $t\mapsto w(t)$ be a smooth curve in $\C^{n-1}$ such that
$w(0)=O$. Let $g_t(\zeta):=(i(1+\zeta)/(1-\zeta)+i\|w(t)\|^2, w(t))$
for $\zeta\in\D$ and $t$ close to $0$. By definition,  $\{g_t\}$ is
a family of complex geodesics in $\Hr^n$, and thus
$\v_t:=G(g_t(\zeta))$ is a smooth real one-parameter family
$\{\v_t\}$ of complex geodesics in $D$ such that
$\v_0(\zeta)=(\zeta,O)$. The associated Jacobi vector field
$J(\zeta)=\frac{\de \v_t}{\de t}(\zeta)$ can be written in the form
\[
J(\zeta)=J_1(\zeta)\frac{\de}{\de z_1}+J^\perp(\zeta),
\]
where $J^\perp(\zeta)=\sum_{k=2}^n J_k(\zeta)\frac{\de}{\de z_k}$
and, since $\v_t(1)=p$ for all $t$, by Corollary \ref{appende} it
follows that $J(1)=O$. Therefore, from \cite[Section 3]{S-T} it
follows that there exist $a\in \C$,  $X, Y\in \C^{n-1}$ (depending
on $J$) and a unique continuous map $M:\oD\to\hbox{GL}(2n-2,\C)$
holomorphic in $\D$  which depends only on $D$ and $\v$ with the
following properties. If $M(\zeta)=\left(\begin{smallmatrix}
M_1(\zeta) & M_2(\zeta)\\M_3(\zeta) & M_4(\zeta)
\end{smallmatrix}\right)$ where the $M_j$'s are suitable $(n-1)\times (n-1)$-matrices
with $M_1(1)=\frac{1}{2}{\sf Id}$, $M_2(1)=\frac{-i}{2}{\sf Id}$
(and $M_3(1), M_4(1)$ satisfy suitable conditions that we do not
need here), then
\begin{equation}\label{J-fatto}
\begin{split}
J_1(\zeta)&=(1-\zeta)(a+\overline{a}\zeta),\\
J^\perp(\zeta)&=i(1-\zeta)(M_1(\zeta)X+M_2(\zeta)Y).
\end{split}
\end{equation}
By the very definition of $G$ and by \eqref{J-fatto}, taking into
account that $G$ maps complex tangent spaces to the boundary of
horospheres in $\Hr^n$ to complex tangent spaces to the boundary of
horospheres in $D$ (see the proof of {\sl Theorem 6.3} in
\cite{B-P}) it follows that for $\Im \xi>0$ and $\zeta\in \D$
\begin{equation}\label{stima-G}
\begin{split}
&\frac{\de G_j}{\de \overline{\xi}}(\xi,0)=0\quad\hbox{for \
}j=1,\ldots, n\\
&\frac{\de G_1}{\de \xi}(\xi,0)=\frac{\de }{\de \xi}(\xi-i)(\xi+i)^{-1}\\
&\frac{\de G_j}{\de \xi}(\xi,0)=0\quad\hbox{for \
}j=2,\ldots, n\\
& \frac{\de G_1}{\de w_j}(\xi,0)=\frac{\de G_1}{\de
\overline{w_j}}(\xi,0)=0\quad\hbox{for\ }j=1,\ldots, n-1.\\
& \frac{\de G_k}{\de
w_j}(i\frac{1+\zeta}{1-\zeta},0)=i(1-\zeta)(M_1(\zeta)S_j+M_2(\zeta)T_j)_k\quad \hbox{for\ }j,k=2,\ldots, n\\
& \frac{\de G_k}{\de
\overline{w}_j}(i\frac{1+\zeta}{1-\zeta},0)=i(1-\zeta)(M_1(\zeta)\overline{S}_j+M_2(\zeta)\overline{T}_j)_k
\quad \hbox{for\ }j,k=2,\ldots, n,\\
\end{split}
\end{equation}
for some vectors $(S_2,\ldots, S_n), (T_2,\ldots, T_n)\in \C^{n-1}$.
Let $S$ (respectively $T$) be the matrix whose columns are
$S_2,\ldots, S_{n-1}$ (respect. $T_2,\ldots, T_{n-1}$) and set
\[
N=\left(
    \begin{array}{cc}
      S & \overline{S} \\
      T & \overline{T} \\
    \end{array}
  \right).
\]
We claim that $N$ is invertible. Indeed, since $dG$ is invertible at
$(i,O)$, equations \eqref{stima-G} imply that the only vector $v$
satisfying $(M_1(0)\ M_2(0))(2\Re \left(\begin{smallmatrix} S\\ T
\end{smallmatrix}\right)v)=0$ is the zero vector $v=O$. Therefore
$S_2,\ldots, S_n, T_2, \ldots, T_n$ form a real basis of $\C^{n-1}$.
From this it follows easily that if the vector $\left(\begin{smallmatrix} v\\
w
\end{smallmatrix}\right)$ belongs to the kernel of $N^t$ then
$v=w=0$ and thus $N$ is  invertible.

Now we are in the good shape to compute $\frac{\de U}{\de
w_j}(\xi,O)$. Since $U=u\circ G=\tl(w)(\|w\|^2-\Im \xi)$, from
\eqref{stima-G} we have for $j=1,\ldots, n-1$ and $\Im \xi>0$
\begin{equation}
\label{scrivo}  -\frac{\de \tl}{\de w_j}(O)\Im
\xi=\sum_{k=2}^{n}[\frac{\de u}{\de z_k}(\frac{\xi-i}{\xi+i},
O)\frac{\de G_k}{\de w_j}(\xi, O)+\overline{\frac{\de u}{\de
z_k}(\frac{\xi-i}{\xi+i}, O)\frac{\de G_k}{\de \overline{w_j}}(\xi,
O)}].
\end{equation}
Notice that since $u$ is plurisubharmonic in $D$ and harmonic on the
complex geodesics whose closure contains $p$, it follows that the
functions $\frac{\de u}{\de z_k}(\frac{\xi-i}{\xi+i}, O)$ are
holomorphic for $\Im\xi>0$. Moreover, by \eqref{stima-G} both
$\frac{\de G}{\de w_j}(\xi,O)$ and $\frac{\de G}{\de
\overline{w}_j}(\xi,O)$ are holomorphic for $\Im\xi>0$. Taking the
real and imaginary part in  \eqref{scrivo} and writing $V(=\frac{\de
\tl}{\de w}(i,O))=C+iD$ with $C,D\in \R^{n-1}$, we find that there
exist two vectors $C', D'\in \R^{n-1}$ such that for all $\Im \xi>0$
\begin{align}
iC_j \xi+iC'_j&=\sum_{k=2}^{n}[\frac{\de u}{\de
z_k}(\frac{\xi-i}{\xi+i}, O)\frac{\de G_k}{\de w_j}(\xi,
O)+\frac{\de u}{\de z_k}(\frac{\xi-i}{\xi+i}, O)\frac{\de
G_k}{\de \overline{w_j}}(\xi, O)],\label{unoG}\\
-D_j \xi+D'_j&=\sum_{k=2}^{n}[\frac{\de u}{\de
z_k}(\frac{\xi-i}{\xi+i}, O)\frac{\de G_k}{\de w_j}(\xi,
O)-\frac{\de u}{\de z_k}(\frac{\xi-i}{\xi+i}, O)\frac{\de G_k}{\de
\overline{w_j}}(\xi, O)]\label{dueG}.
\end{align}
Let $V'=iC'+D'$, let $f_k(\zeta)=-2i(1-\zeta)^2\frac{\de u}{\de
z_k}(\zeta,O)$ and let $f=(f_1,\ldots, f_n)$ for $\zeta\in\D$.
Summing (respectively subtracting) \eqref{unoG} with \eqref{dueG},
composing with $\zeta\mapsto i\frac{1+\zeta}{1-\zeta}$, multiplying
by $(1-\zeta)$ and using \eqref{stima-G} we obtain for
$\zeta\in\oD\setminus\{1\}$
\[
\left(
  \begin{array}{c}
    \zeta(V+V')+(V-V')  \\
    \zeta(\overline{V}-\overline{V}')+(\overline{V}+\overline{V}')\\
  \end{array}
\right)=N^t \left(
              \begin{array}{c}
                M_1(\zeta)^t \\
                M_2(\zeta)^t \\
              \end{array}
            \right)\cdot f(\zeta).
\]
From this, since $N$ is invertible and also $M_1(\zeta), M_2(\zeta)$
are invertible for $\zeta$ close to $1$ (since by the very
definition $M_1(1)=\frac{1}{2}{\sf Id}$ and $M_2(1)=\frac{-i}{2}{\sf
Id}$) it follows that $f(\zeta)$ has a limit $L$ at $\zeta=1$ and
\begin{equation}\label{infine}
\left(
  \begin{array}{c}
    4V \\
    4\overline{V} \\
  \end{array}
\right)=N^t\left(
  \begin{array}{c}
    {\sf Id} \\
    -i{\sf Id} \\
  \end{array}
\right)L.
\end{equation}
Therefore $(S^t-iT^t)L-(S^t+iT^t)\overline{L}=O$. Writing
$L=\alpha+i\beta$ for $\alpha,\beta\in \R^{n-1}$, we have
$S^t\beta-T^t\alpha=O$ and, since $\alpha, \beta$ are real, this is
equivalent to
\[
N^t\left(
     \begin{array}{c}
       \beta  \\ -\alpha
     \end{array}
   \right)=O.
\]
But $N$ is invertible and therefore $\alpha=\beta=O$ which means
$L=O$. Finally, from \eqref{infine} it follows that $V=O$.
\end{proof}

The pluricomplex Poisson kernel can be also characterized in terms
of its level sets:

\begin{prop}\label{level-unico}
Let $D\subset \C^n$ be a bounded strongly convex domain with smooth
boundary and let $p\in \de D$. Let $u\in\ps(D)\cap
L_{loc}^\infty(D)$ be such $(\de\debar u)^n=0$ in $D$ and
$\lim_{z\to x}u(z)=0$ for all $x\in \de D\setminus\{p\}$. If $u$ has
the same level sets of $\Omega_{D,p}$ then there exists $c>0$ such
that $u=c\Omega_{D,p}$.
\end{prop}
\begin{proof}
By hypothesis there exists a function $Y:\R^-\to \R^-$ such that
$u(z)=Y(\Omega_{D,p}(z))$ for all $z\in D$. We need to show that
there exists $c>0$ such that $Y(t)=ct$ for all $t\in\R^-$. To this
aim, since each complex geodesic whose closure contains $p$
intersects every horosphere, it is enough to prove that
$u(z)=c\Omega_{D,p}(z)$ for $z$ belonging to any complex geodesic
whose closure contains~$p$.

Let $S$ be a complex geodesic in $D$ such that $p\in \overline{S}$
and $\rho:D\to S$ the associated Lempert's projection. We can
 assume that $D$ is linearizated along $S$ in Lempert's
coordinates. Let $\tilde{B}$ be a open disc relatively compact in
$S$. Let
\[
\mathcal P:=\begin{cases} \tilde{v}\in \subh (\tilde{B})\\
\limsup_{\zeta\to x}\tilde{v}(\zeta)\leq u(x)\quad \forall x\in \de
\tilde{B}
\end{cases}
\]
If we prove that $u|_{\tilde{B}}$ is the maximum of $\mathcal P$
then, by the arbitrariness of $\tilde{B}$ it follows that $u$ is
harmonic on $S$. Therefore $u\circ \v$ is harmonic and negative in
$\D$ and it is zero on $\de \D\setminus\{1\}$, hence it is a
constant multiple of the Poisson kernel of $\D$. That is, there
exists $c>0$ such that $u(\v(\zeta))=c\Omega_{D,p}(\v(\zeta))$ for
all $\zeta\in\D$, as wanted.

In order to prove that $u|_{\tilde{B}}$ is the maximum of $\mathcal
P$, let $\epsilon>0$ small. Let $T=\rho^{-1}(\tilde{B})\cap D$ and
let $B=\{z\in T: \hbox{dist}(z, \de D)>\epsilon\}$ (a cylinder in
$D$). The boundary of the set  $B$ is made of two parts: $R_1$ which
has the property that $\rho(R_1)=\de \tilde{B}$ and $R_2$ (the
bottom and top of the cylinder) such that $\rho(R_2)\subset
\tilde{B}$; $\de B=R_1\cup R_2$. Since $u=0$ on $\de D$ and
$p\not\in T$, then we can choose $\epsilon$ so small that
$\inf_{x\in R_2}u(x)>\max_{x\in \de \tilde{B}}u(x)$.

Let $\tilde{v}\in \mathcal P$. Let $v:=\tilde{v}\circ \rho|_{B}$.
Then $v$ is plurisubharmonic in $B$ and $\sup_{x\in
B}v(x)=\sup_{x\in \de \tilde{B}}(\limsup_{z\to x}v(z))$. In
particular by construction $\limsup_{z\to x} v(z)\leq u(x)$ for all
$x\in R_2$. Also, we have that $\limsup_{B\ni z\to
x}v(z)=\limsup_{B\ni z\to x}\tilde{v}(\rho(z))\leq u(\rho(x))$ for
all $x\in R_1$. Now $u$ has the same level sets of $\Omega_{D,p}$
and thus by \eqref{ins2} we have that $u(x)\geq u(\rho(x))$ for all
$x\in D$ and hence $\limsup_{B\ni z\to x}v(z)\leq u(x)$ for all
$x\in R_1$. Therefore $\limsup_{B\ni z\to x}v(z)\leq u(x)$ for all
$x\in \de B$ and by the maximality of Monge-Ampere solutions, it
follows that $v\leq u$ in $B$ and in particular $\tilde{v}\leq
u|_{\tilde{B}}$ and the arbitrariness of $\tilde{v}$ implies that
$u|_{\tilde{B}}$ is maximal in $\mathcal P$.
\end{proof}

The  previous argument, together with Theorem
\ref{foliazione-unica}, shows that if $u\in\ps(D)\cap C^2(D)$ is
such that $(\de\debar u)^n=0$ on $D$ and $\lim_{z\to x}u(z)=0$ for
all $x\in \de D\setminus\{p\}$ then $u=c\Omega_{D,p}$ for some $c>0$
if and only if $u(\rho(z))\leq u(z)$ for all $z\in D$ and for all
Lempert's projections $\rho$.

More generally, if $f:D\to D$ is holomorphic and $f(p)=p$ as
non-tangential limit we can define the {\sl boundary dilatation
coefficient} $\alpha_f(p)$ of $f$ at $p$ by means of
\[
\frac{1}{2}\log \alpha_f(p):=\liminf_{z\to p}[k_D(z_0,
z)-k_D(z_0,f(z))].
\]
It turns out that $\alpha_f(p)>0$  and, if $\alpha_f(p)<\infty$, we
can rephrase  Abate's generalization of the classical Julia Lemma
(see \cite{Ab}, \cite{Ab2}) saying that
$\alpha_f(p)f^\ast(\Omega_{D,p})\leq \Omega_{D,p}$. In \cite[Theorem
7.3]{B-P}, with a slightly more regularity assumption required on
$f$ at $p$, it is proved that $f$ is an automorphism of $D$ if and
only if $f^\ast(\Omega_{D,p})= \alpha_f(p)\Omega_{D,p}$.

Using Abate's version of the Julia-Wolff-Caratheodory theorem for
strongly convex domains (see \cite{Ab}, \cite{Ab3}) it is easy to
see that  $\alpha_\rho(p)= 1$ for all Lempert's projections $\rho$.
Therefore,  the above discussion shows that the property
$f^\ast(\Omega_{D,p})\leq \alpha_f(p)\Omega_{D,p}$ characterizes
$\Omega_{D,p}$. In other words:

\begin{prop}\label{contra-unico}
Let $D\subset \C^n$ be a bounded strongly convex domain with smooth
boundary and let $p\in \de D$. Let $u\in\ps(D)\cap C^2(D)$ be such
that $(\de\debar u)^n=0$ in $D$ and $\lim_{z\to x}u(z)=0$ for all
$x\in \de D\setminus\{p\}$. Then there exists $c\geq 0$ such that
$u=c\Omega_{D,p}$ if and only if for all $f:D\to D$ holomorphic such
that $f(p)=p$ as non-tangential limit and $\alpha_f(p)<\infty$ it
follows that
\[
\alpha_f(p)f^\ast(u)\leq  u.
\]
\end{prop}

Some remarks about uniqueness properties are in order. First, it
would be interesting to see whether Theorem \ref{foliazione-unica}
(and thus its corollaries) holds without any regularity hypothesis
on $u$. A direct argument using the sub-media property of
plurisubharmonic functions shows that Theorem \ref{foliazione-unica}
holds in the unit ball $\B^n$ with no regularity hypothesis on $u$.
Such an argument seems however to fail in general.

Another (maybe more) interesting open question is the following:

\begin{question}Let $D\subset \C^n$ be a bounded strongly convex domain with
smooth boundary and let $p\in \de D$. Let $u\in \ps(D)\cap
L_{loc}^\infty(D)$ be such that $(\de \debar u)^n=0$ in $D$ and
$\lim_{z\to x}u(z)=0$ for all $x\in \de D\setminus\{p\}$. Is it true
that $u=c\Omega_{D,p}$ for some constant $c\geq 0$?
\end{question}

As we already recalled, the answer to such a question is ``yes'' in
case $D=\D$ the unit disc, $u<0$ in $\D$ and $\Omega_{D,p}$ is the
(negative) Poisson kernel.

\section{Reproducing formulas}\label{riproduce-sezione}

Let $D$ be a bounded strongly convex domain in $\C^n$ with smooth
boundary. As usual, let $d^c:=i(\debar-\de)$. Let $r$ be a defining
function of $D$ and let $\omega_D$ be the real $(2n-1)$-form defined
as
\[
\omega_{\de D}:=\frac{(dd^c r)^{n-1}\wedge d^c r}{\|dr\|^n}|_{\de
D}.
\]
such a form $\omega_{\de D}$ is positive  and it is easily seen to
be independent of the defining function $r$ chosen to define it.

Let $L_{D,z_0}$ denote the Lempert solution of
\eqref{monge-inner-intro} and denote by $\Omega_{D,p}$  the solution
of \eqref{monge-boundary} with singularity at $p\in \de D$ given by
Theorem \ref{filo-giorgio}. From the very definition of
$\Omega_{D,p}$  and since the boundary spherical representation
$\Phi_p$ of Chang-Hu-Lee is smooth out of the diagonal of $\de
D\times \de D$ as the vertex $p$ varies on $\de D$ (see
\cite[Theorem 3]{CHL}) it follows that the map $\overline{D}\times
\de D\ni (z,p)\mapsto \Omega_{D,p}(z)\in \R$ is $C^\infty$ on
$(\overline{D}\times \de D)\setminus \{(p,p)\in \de D\times \de
D\}$.

We briefly recall Demailly's theory \cite{De1}, \cite{De}. Let
$\v\in \ps(D)$ be such that $\exp(\v)\in C^0(\overline{D})$, that
$\v<0$ on $D$ and that $\v=0$ on $\de D$. Let $R<0$ and let
$B_R=\{z\in D: \v(z)<R\}$. Moreover let $S_R=\de B_R$ and
$\v_R(z)=\max\{\v(z), R\}$. By \cite[(1.4)]{De} we can write
\[
(dd^c \v_R)^n={\bf 1}_{\C^n\setminus B_R}(dd^c \v)^n+\mu_{\v,R}
\]
where ${\bf 1}_{\C^n\setminus B_R}$ is the characteristic function
of $\C^n\setminus B_R$ and $\mu_{\v,R}$ is a positive measure
supported on $S_R$. By \cite[Th\'eor\`eme 3.1]{De} if the total
Monge-Amp\`ere mass of $\v$ is finite, {\sl i.e.}, if $\int_D (dd^c
\v)^n<+\infty$, then as $R\to 0$ the measures $\mu_{\v,R}$ converge
weakly on $\C^n$ to a positive measure $\mu_\v$ supported on $\de
D$, with total mass $\int_D (dd^c \v)^n$. We denote by $\mu_z$ the
limit measure of $L_{D,z}$. By \cite[Th\'eor\`eme 5.1]{De} it
follows that for all $F\in\ps(D)\cap C^0(\overline{D})$ we have the
following representation formula:
\begin{equation}\label{demail}
F(z)=\mu_z(F)-\frac{1}{2\pi^n}\int_{w\in D} |L_{D,z}(w)|\; dd^c
F(w)\wedge (dd^c L_{D,z})^{n-1}(w).
\end{equation}

We can prove the following result:
\begin{teo}\label{comemu}
Let $D$ be a bounded strongly convex domain in $\C^n$ with smooth
boundary. Then
\[
d\mu_z(p)=|\Omega_{D,p}(z)|^n  \omega_{\de D}(p).
\]
\end{teo}
\begin{proof}
First of all, since $L_{D,z}$ is $C^\infty$ on
$\overline{D}\setminus\{z\}$ and $dL_{D,z}|_{\de D}\neq 0$, arguing
as in \cite{De} we see that
\begin{equation*}
d\mu_z=(dd^c L_{D,z})^{n-1}\wedge d^c L_{D,z}|_{\de D}.
\end{equation*}
From \eqref{legati} we have
\[
|\Omega_{D,p}(z)|=\|\frac{\de L_{D,z}}{\de \nu_p}(p)\|=\|d
(L_{D,z})_p\|,
\]
where the last equality follows from $L_{D,z}|_{\de D}= 0$ which
implies that $d (L_{D,z})_p$ is a positive multiple of $\nu_p$, the
unit normal to $\de D$ at $p\in \de D$ (here, as usual and with an
abuse of notation, we identify the gradient of a function with its
differential). Thus
\[
d\mu_z=|\Omega_{D,p}(z)|^n\frac{(dd^c L_{D,z})^{n-1}\wedge d^c
L_{D,z}}{\|dL_{D,z}\|^n}|_{\de D}.
\]
To end  the proof we only need to check that
\[
\omega_{\de D}=\frac{(dd^c L_{D,z})^{n-1}\wedge d^c
L_{D,z}}{\|dL_{D,z}\|^n}|_{\de D}.
\]
To this aim, it is enough to show that if $r$ is a (local) defining
function for $D$ on a neighborhood $U_p$ of  $p\in \de D$, then
$L_{D,z}=h\cdot r$ on $U_p\cap \overline{D}$ for some positive $h\in
C^\infty(U_p\cap \overline{D})$ for then a direct computation gives
the result. Up to changes of coordinates we can assume that $U_p\cap
D=\{(x,y)\in \C\times \C^{n-1}: x<0\}$. Thus $L_{D,z}(x,y)/x$ is
defined and positive on $U_p\cap D$. If we let $h(x,y)=\int_0^1
\frac{\de L_{D,z}}{\de x}(tx,y)dt$ then $h$ is $C^\infty(U_p\cap
\overline{D})$ and coincides with $L_{D,z}(x,y)/x$ in $U_p\cap D$.
Moreover, since $dL_{D,z}\neq 0$ on $\de D$ it follows that $h>0$ on
$U_p\cap \overline{D}$.
\end{proof}

From \eqref{demail} and Theorem \ref{comemu} we obtain:

\begin{teo}\label{riproduce}
Let $F\in \ps(D)\cap C^0(\overline{D})$. Then for all $z\in D$
\begin{equation*}
\begin{split}
F(z)&=\int_{p\in\de D} |\Omega_{D,p}(z)|^n F(p) \omega_{\de
D}(p)\\&-\frac{1}{2\pi^n}\int_{w\in D} |L_{D,z}(w)|\; dd^c
F(w)\wedge (dd^c L_{D,z})^{n-1}(w).
\end{split}
\end{equation*}
In particular if $F$ is pluriharmonic then
\[
F(z)=\int_{p\in\de D} |\Omega_{D,p}(z)|^n F(p) \omega_{\de D}(p).
\]
\end{teo}

\begin{nota}
If $F\in  C^2(\overline{D})$ (but {\sl not} plurisubharmonic in $D$)
then there exists $C>0$ such that $F(z)+C\|z\|^2\in \ps(D)\cap
C^0(\overline{D})$. Thus Theorem \ref{riproduce} applies and one
gets
\begin{equation*}
\begin{split}
F(z)+C\|z\|^2&=\int_{p\in\de D} |\Omega_{D,p}(z)|^n F(p) \omega_{\de
D}(p)+C\int_{p\in\de D} |\Omega_{D,p}(z)|^n \|p\|^2 \omega_{\de
D}(p)\\&-\frac{1}{2\pi^n}\int_{w\in D} |L_{D,z}(w)|\; dd^c
F(w)\wedge (dd^c L_{D,z})^{n-1}(w)\\&-C\frac{1}{2\pi^n}\int_{w\in D}
|L_{D,z}(w)|\; dd^c \|w\|^2\wedge (dd^c
L_{D,z})^{n-1}(w)\\&=\int_{p\in\de D} |\Omega_{D,p}(z)|^n F(p)
\omega_{\de D}(p)\\&-\frac{1}{2\pi^n}\int_{w\in D} |L_{D,z}(w)|\;
dd^c F(w)\wedge (dd^c L_{D,z})^{n-1}(w)+C\|z\|^2.
\end{split}
\end{equation*}
Therefore Theorem \ref{riproduce} applies to any $F\in
C^2(\overline{D})$ (not necessarily plurisubharmonic). As a
consequence it follows that the kernel $|\Omega_{D,p}(z)|^n
\omega_{\de D}(p)$ is the {\sl unique} reproducing kernel associated
to $L_{D,z}$, namely,
 \eqref{demail} cannot hold  with any other measure $T_z$ in place of $\mu_z$.
\end{nota}

\end{document}